\documentclass[12pt]{amsart}
\usepackage{amssymb}
\usepackage{amsthm}
\usepackage{amsmath}
\usepackage{enumerate}
\usepackage{cancel}
\usepackage[normalem]{ulem}
\usepackage[utf8]{inputenc}

\newtheorem{theorem}{Theorem}

\newtheorem{lemma}{Lemma}
\newtheorem{proposition}{Proposition}

\newtheorem{corollary}{Corollary}

\begin{document}
	\title[Topological loops having decomposable solvable multiplication group]
	{Topological loops having decomposable solvable multiplication group}
	
	\thanks{The paper is supported by the National Research, Development and Innovation Office (NKFIH) Grant No. K132951 and by the EFOP-3.6.1-16-2016-00022 project. The later project is co-financed by the European Union and the European Social Fund. }
	
	\author[A. Al-Abayechi]{Ameer Al-Abayechi} \address{Institute of Mathematics, University of Debrecen, and Doctoral School of Mathematical and Computational Sciences of University of Debrecen, Hungary} \email{ameer@science.unideb.hu}

	\author[\'A. Figula]{\'Agota Figula}
	\address{Institute of Mathematics, University of Debrecen, Hungary} \email{figula@science.unideb.hu}
	\email{ORCID: 0000-0002-8095-6074}

	\subjclass[2010]{20N05, 22E25, 17B30, 57M60, 57S20}
	
	\keywords{topological loop, multiplication group and inner mapping group, solvable Lie groups and Lie algebras, connected transversals, centrally nilpotent loops of class two}
	
	\date{}

\begin{abstract} 
In this paper we deal with the class $\mathcal{C}$ of decomposable solvable Lie groups having
dimension at most six. We determine those Lie groups in $\mathcal{C}$
and their subgroups which are the multiplication group $Mult(L)$ and the inner
mapping group $Inn(L)$ for three-dimensional connected simply connected
topological loops $L$. These loops $L$ have one- or two-dimensional centre and their group $Mult(L)$ has two- or
three-dimensional commutator subgroup. Together
with this result we obtain that every at most $3$-dimensional connected
topological proper loop having a solvable Lie group of dimension at most six as
its multiplication group is centrally nilpotent of class two.
\end{abstract}



\maketitle

{\small\centerline{In memoriam Prof. Dr. Heinrich Wefelscheid (16.4.1941-18.4.2020)}}

\bigskip

\section{Introduction}
\label{introduction}
The notion of the multiplication group $Mult(L)$ and the inner mapping group $Inn(L)$ of a loop $L$ was introduced and firstly investigated by A. A. Albert and R. H. Bruck. Since their papers \cite{albert1}-\cite{bruck} much work has been done to study the correspondences between the structure of the loop $L$ and that of the group $Mult(L)$ and $Inn(L)$. In particular many results relate nilpotence and solvability of loops to the analogous properties of their groups $Mult(L)$ and
$Inn(L)$ (\cite{mazur}, \cite{niemenmaa}, \cite{vojtechovsky}, \cite{vesanen}).
M. Niemenmaa and T. Kepka established in \cite{kepka} the necessary and sufficient conditions for a group $G$ to be the multiplication group of $L$. In their criteria the existence of special transversals $A$ and $B$ with respect to a subgroup $K$ of $G$ plays an important role. These transversals belong to the sets of left and right translations of $L$ whereas
$K$ corresponds to the inner mapping group of $L$. For finite loops the importance of
the permutation groups $Mult(L)$ and $Inn(L)$ as well as the connected transversals A and B is documented in many papers (cf.
\cite{niem4}-\cite{drapal}, \cite{nagyg}, \cite{niem0}, \cite{niem1}, \cite{vesanen1}).

Topological and differentiable loops are investigated thoroughly by P. T. Nagy and K. Strambach in \cite{loops} as continuous and differentiable sections in Lie groups. Part II of \cite{loops} is devoted for explicit description and determination of topological and smooth loops on low dimensional manifolds.
Following their approach this article is a contribution to the classification of connected topological loops $L$ of dimension $3$ having a solvable Lie group as their multiplication group. Each $2$-dimensional connected topological proper loop having a Lie group as its multiplication group has nilpotency class two (cf. \cite{figula0}). This nilpotency property is valid for $3$-dimensional connected topological loops $L$ having either a solvable Lie group of dimension at most $5$ or a $6$-dimensional indecomposable solvable Lie group as their group $Mult(L)$ (cf. \cite{figula1}-\cite{figula2}, \cite{figula4}-\cite{figula3}).
Furthermore, in the class of the at most $5$-dimensional solvable non-nilpotent Lie groups only decomposable groups occur as the multiplication group of $L$ (cf. \cite{figula1}).
In this paper we show that the centrally nilpotency of class two property is satisfied for $3$-dimensional topological loops $L$ if the group $Mult(L)$ is a $6$-dimensional decomposable solvable Lie group (cf.  Propositions \ref{1dimcentre},
\ref{2dimensionalcentre}).

The $3$-dimensional connected simply connected topological loops $L$ having a solvable indecomposable Lie group of dimension at most $6$ as their multiplication group are known (cf. \cite{figula2}-\cite{ffaa}, \cite{figula3}).
To complete this classification for every solvable Lie group of dimension at most $6$ we determine
the multiplication Lie groups and the inner mapping subgroups for $L$ in the class of the decomposable solvable Lie groups (cf. Theorems \ref{multiplication1}, \ref{multiplication2}, Proposition \ref{nilpotent}).

After the presentation of the necessary concepts, Proposition
\ref{6dimdecomtrialcentre} shows that the decomposable solvable Lie groups of dimension
$6$ with discrete centre are not the multiplication group of $3$-dimensional connected topological loops. In section $4$, respectively
in $5$ we treat the decomposable solvable Lie groups having a $1$-dimensional, respectively a $2$-dimensional centre. Among the
$6$-dimensional decomposable solvable Lie groups with $1$-dimensional centre there are $18$ families of Lie groups which are multiplication groups of $L$. In the class of the $6$-dimensional decomposable solvable Lie groups with $2$-dimensional centre $9$ families can be represented as the group $Mult(L)$ of $L$. All these Lie groups have $3$-dimensional commutator subgroup (see Corollary \ref{cor1}) and depend on at most two real parameters.

\section{Preliminaries}
\label{sec:1}
A set $L$ with a binary operation $(x,y) \mapsto x \cdot y$ is called a loop if there exists an element
$e \in L$ such that $x=e \cdot x=x \cdot e$ holds for all $x \in L$ and for each $x \in L$ the left translations
$\lambda _x:L \to L, \lambda_x (y)=x \cdot y$ and the right translations $\rho _x:L \to L$, $\rho_x(y)=y \cdot x$
 are bijections of $L$. A loop $L$ is proper if it is not a group. The maps
$(x,y) \mapsto x \backslash y=\lambda_x^{-1}(y)$, respectively $(x,y) \mapsto y /x= \rho_x^{-1}(y)$, $x,y \in L$ are further binary operations on $L$. The permutation group $Mult(L)=\langle \lambda_x, \rho_x; x \in L \rangle $ is called the multiplication group of $L$. The inner mapping group $Inn(L)$ of $L$ is the stabilizer of the identity element $e \in L$ in $Mult(L)$.

Let $G$ be a group,
let $K \le G$, and let $A$ and $B$ be two left transversals to $K$ in $G$. We say
that $A$ and $B$ are $K$-connected if $a^{-1}b^{-1}ab \in K$ for every $a \in A$ and $b \in B$.
The core $Co_G(K)$ of $K$ in $G$ is the largest normal subgroup of $G$ contained
in $K$. If $L$ is a loop, then $\Lambda(L) = \{\lambda_x; x \in L\}$ and $R(L) =\{\rho_x; x \in L\}$
are $Inn(L)$-connected transversals in the group $Mult(L)$.
Theorem 4.1 in \cite{kepka} states the following necessary and sufficient conditions for a group $G$ to be the multiplication group of a loop $L$:

\begin{proposition} \label{kepka}
A group $G$ is isomorphic to the multiplication group of a loop
if and only if there exists a subgroup $K$ with $Co_G(K) = 1$ and $K$-connected left
transversals $A$ and $B$ satisfying $G =\langle A, B \rangle $.
\end{proposition}

The kernel of a homomorphism $\alpha :(L, \cdot ) \to (L', \ast )$ of a loop $L$ into a loop $L'$ is a normal subloop $N$ of $L$.
The centre $Z(L)$ of a loop $L$ consists of all elements $z$ which satisfy the equations
$z x \cdot y=z \cdot x y, \ x \cdot y z=x y \cdot z, \ x z \cdot y=x \cdot z y, \ z x =x z$ for  all $x,y \in L$.
If we put $Z_0=e$, $Z_1=Z(L)$ and $Z_i/Z_{i-1}=Z(L/Z_{i-1})$, then we obtain a series of normal subloops of $L$.
If $Z_{n-1}$ is a proper subloop of $L$ but $Z_n=L$, then $L$ is centrally nilpotent of class $n$.

A loop $L$ is called classically solvable if there is a series $\{e\}=L_0 \le L_1 \le \dots \le L_n=L$ of subloops of $L$ such that, for every $i=1,2, \cdots ,n$, $L_{i-1}$ is normal in $L_{i}$ and each factor loop $L_{i}/L_{i-1}$ is a commutative group.

We often use the following Lemma which is proved in Theorems 3, 4 and 5 of \cite{albert1}, in Lemma 1.3, IV.1 of \cite{bruck}, in Lemma 2.3 of \cite{figula2} and in Proposition 2.7 of \cite{kepka}.

\begin{lemma}\label{bruck} Let $L$ be a loop with multiplication group $Mult(L)$, inner mapping group $Inn(L)$ and identity element
$e$.
\begin{enumerate}
  \item Let $\alpha $ be a homomorphism of the loop $L$ onto the loop $\alpha (L)$ with kernel $N$. Then $\alpha $ induces a homomorphism of the group $Mult(L)$ onto the group $Mult(\alpha (L))$.
Let $M(N)$ be the set $\{ m \in Mult(L); \ x N=m(x) N \ \hbox{for  all} \ x \in L \}$.
Then $M(N)$ is a normal subgroup of $Mult(L)$ containing the multiplication group $Mult(N)$ of the loop $N$ and the multiplication group of the factor loop $L/N$ is isomorphic to $Mult(L)/M(N)$.
  \item For every normal subgroup $\mathcal{N}$ of $Mult(L)$ the orbit $\mathcal{N}(e)$ is a normal subloop of $L$ and
$\mathcal{N} \le  M(\mathcal{N}(e))$.
  \item The core of $Inn(L)$ in $Mult(L)$ is trivial and the normalizer of $Inn(L)$ in $Mult(L)$ is the direct product
$Inn(L) \times Z$, where $Z$ is the centre of the group $Mult(L)$.
\end{enumerate}
\end{lemma}

A loop $L$ is called topological  if $L$ is a topological space and the binary operations
$(x,y) \mapsto x \cdot y, \ (x,y) \mapsto x \backslash y, (x,y) \mapsto y/x :L \times L \to L$  are continuous.
In general the multiplication group of $L$ is a topological transformation group that does not have a natural (finite dimensional) differentiable structure. We investigate a $3$-dimensional connected topological loop having a solvable Lie group as its multiplication group. The first assertion of the following Lemma is proved in \cite{hofmann2}, IX.1, the second assertion is showed in \cite{figula1}, Lemma 5.

\begin{lemma} \label{simplyconnected} For every connected topological loop there exists the universal covering loop
$L$. If $L$ is a $3$-dimensional connected simply connected topological loop having a solvable Lie group as its multiplication group, then it is homeomorphic to $\mathbb R^3$.
\end{lemma}

The elementary filiform Lie group ${\mathcal F}_{n}$ is the simply connected nilpotent Lie group of dimension $n \ge 3$ whose Lie algebra is elementary filiform, i.e. it has a basis $\{ e_1, \cdots , e_{n} \}$  with $[e_1, e_i]= e_{i+1}$ for
$2 \le i \le n-1$.  A $2$-dimensional simply connected topological loop $L_{\mathcal F}$ is called an elementary filiform loop if its multiplication group is an elementary filiform group ${\mathcal F}_{n}$ of dimension $n \ge 4$ (\cite{figula6}).

A Lie algebra is called decomposable, if it is the direct sum of two proper ideals.
In this paper we assume that the multiplication group of $L$ is a $6$-dimensional solvable decomposable Lie group or a nilpotent decomposable Lie group of dimension $\le 5$. The next lemma summarizes the known results about the $3$-dimensional topological loops having solvable decomposable Lie groups as their multiplication groups (cf. Lemmata 3.4, 3.5, 3.6 and Propositions 3.7, 3.8 in \cite{figula6}, pp. 390-393, Theorem 11 in \cite{albert1}, Theorem 6, Sections 4 and 5 in \cite{figula1}, Propositions 2.6, 2.7 in \cite{figula2}, Lemma 6 (d) in \cite{figula3}, Chapter I in \cite{bruck}) which are often used in the paper.

\begin{lemma} \label{lemmaregi}
Let $L$ be a $3$-dimensional proper connected simply connected topological loop such that its multiplication group
$Mult(L)$ is a $6$-dimensional solvable decomposable Lie group.
\newline
\noindent
a) Then the centre $Z$ of the group $Mult(L)$ and the centre $Z(L)=Z(e)$ of the loop $L$, where $e$ is the identity element of
 $L$, are isomorphic. Moreover, the centre $Z$ has dimension $\le 2$.
\newline
\noindent
b) The loop $L$ is classically solvable and it has a $1$-dimensional connected normal subloop $N$. Every such subloop $N$ of $L$ is
isomorphic to $\mathbb R$ and lies in a $2$-dimensional connected normal subloop $M$ of $L$. The factor loop $L/M$ is isomorphic to $\mathbb R$, whereas the loop $M$ and the factor loop $L/N$ are isomorphic either to the Lie group $\mathbb R^2$ or to the
$2$-dimensional non-abelian Lie group ${\mathcal L}_2$ or to an elementary filiform loop $L_{\mathcal F}$.
\newline
\noindent
c) If $Mult(L)$ has discrete centre, then for every normal subloop $N \cong \mathbb R$ of $L$ the factor loop $L/N$ is isomorphic either to the group ${\mathcal L}_2$ or to a loop $L_{\mathcal F}$. The group $Mult(L)$ has a normal subgroup $S$ containing
$Mult(N) \cong \mathbb R$ such that the factor group $Mult(L)/S$ is isomorphic to the direct product
${\mathcal L}_2 \times {\mathcal L}_2$ if $L/I(e) \cong {\mathcal L}_2$ or to an elementary filiform Lie group
${\mathcal F}_{n}$, $n=4, 5$, if $L/I(e) \cong L_{\mathcal F}$. The normal subloop $M$ containing $N$ is isomorphic either to
$\mathbb R^2$ or to ${\mathcal L}_2$ or to $L_{\mathcal F}$. The group
$Mult(L)$ has a normal subgroup $V$ such that the orbit $V(e)$ is the loop $M$,
$Mult(L)/V \cong \mathbb R$, $V$ contains the inner mapping group $Inn(L)$ of $L$, the group $Mult(M)$ of $M$ and the commutator subgroup of $Mult(L)$.
\newline
\noindent
d) If $\hbox{dim} (Z(L))=1$, then for every normal subloop $N \cong \mathbb R$ of $L$ we have one of the following possibilities:
\newline
\noindent
(i) The factor loop $L/N$ is isomorphic to $\mathbb R^2$. Then $L$ is centrally nilpotent of class $2$, $N$ coincides with the centre $Z(L)$ of $L$ and the group $Mult(L)$ is a semidirect product of the normal subgroup
$P=Z \times Inn(L) \cong \mathbb R^4$ by a group
$Q \cong \mathbb R^2$ such that the orbit $P(e)$ is $Z(L)$.
\newline
\noindent
(ii) The factor loop $L/N$ is isomorphic either to the Lie group ${\mathcal L}_2$ or to a loop $L_{\mathcal F}$. Then case c) is fulfilled.
In particular, if $N=Z(L)$, then $M$ is not isomorphic to the group ${\mathcal L}_2$.
\newline
\noindent
e) If $\hbox{dim} (Z(L))=2$, then $L$ is centrally nilpotent of class $2$ and the group $Mult(L)$ is a semidirect product of the normal subgroup $V=Z \times Inn(L) \cong \mathbb R^5$ by a group $Q \cong \mathbb R$, where
$\mathbb R^2 =Z \cong Z(L)$. $V$ contains the commutator subgroup $Mult(L)'$ of $Mult(L)$. The group $Mult(L)$ is either nilpotent or it has $5$-dimensional abelian nilradical. For every $1$-dimensional connected subgroup $I$ of $Z$ the orbit $I(e)$ is a connected central subgroup of $L$ isomorphic to $\mathbb R$ and the factor loop $L/I(e)$ is isomorphic either to the group
$\mathbb R^2$ or to an elementary filiform loop $L_{\mathcal F}$. If
$L/I(e) \cong \mathbb R^2$, then case d) (i) holds.
If $L/I(e) \cong L_{\mathcal F}$, then the group $Mult(L)$ has a normal subgroup $S$ containing $I \cong \mathbb R$  such that 
one has $S(e)=I(e)$ and the factor group $Mult(L)/S$ is isomorphic either to
${\mathcal F}_{4}$ or to ${\mathcal F}_{5}$. The case $Mult(L)/S \cong {\mathcal F}_{5}$ occurs only if
$Mult(L)=\mathbb R \times {\mathcal F}_{5}$.
\newline
\noindent
f) If the group $Mult(L)$ is nilpotent, then $L$ is centrally nilpotent.
\end{lemma}

The next Lemma, which is proved in \cite{figula4}, Proposition 3.3, is a useful tool to exclude those Lie algebras which are not the Lie algebra of the multiplication group of a
$3$-dimensional topological loop.

\begin{lemma} \label{tool}
Let $L$ be a $3$-dimensional proper connected simply connected topological loop having a $6$-dimensional solvable Lie
algebra ${\bf g}$ as the Lie algebra of its multiplication group.
\newline
\noindent
a) For each $1$-dimensional ideal ${\bf i}$ of ${\bf g}$ the orbit $I(e)$, where $I$ is the simply connected Lie group of 
${\bf i}$ and $e$ is the identity element of $L$, is a normal subgroup of $L$ isomorphic to
$\mathbb R$. We have one of the following possibilities:
\newline
\noindent
(i) The factor loop $L/I(e)$ is isomorphic to $\mathbb R^2$. Then ${\bf g}$ contains the ideal
${\bf p}={\bf z} \oplus {\bf inn(L)} \cong \mathbb R^4$ such that the commutator ideal ${\bf g}'$ of ${\bf g}$ lies in
${\bf p}$. Here ${\bf z}$ is the $1$-dimensional centre of ${\bf g}$ and ${\bf inn(L)}$ is the Lie algebra of the inner mapping group
$Inn(L)$.
\newline
\noindent
(ii) The factor loop $L/I(e)$ is isomorphic either to the group ${\mathcal L}_2$ or to a loop $L_{\mathcal F}$. Then there exists an ideal ${\bf s}$ of ${\bf g}$ such that ${\bf i} \le {\bf s}$ and the factor Lie algebra ${\bf g}/{\bf s}$ is isomorphic either to  ${\bf l}_2 \oplus {\bf l}_2$, where ${\bf l}_2$ is the $2$-dimensional solvable non-abelian Lie algebra, or to an elementary filiform Lie algebra ${\bf f}_n$, $n=4, 5$.
\newline
\noindent
b) Let ${\bf a}$ be an ideal of ${\bf g}$ such that $\hbox{dim}({\bf a})=2$,
${\bf a} \subseteq {\bf g}'$ and
the factor Lie algebra ${\bf g}/{\bf a}$ is isomorphic neither to ${\bf l}_2 \oplus {\bf l}_2$ nor to ${\bf f}_4$. Then the orbit
$A(e)$, where $A$ is the simply connected Lie group of ${\bf a}$, is either a $2$-dimensional connected normal subloop $M$ of $L$ or the factor loop $L/A(e)$ is isomorphic to $\mathbb R^2$.

If one has $A(e)=M$, then there exists a
$5$-dimensional ideal ${\bf v}$ of ${\bf g}$ containing the Lie algebra
${\bf inn(L)}$, the Lie algebra ${\bf mult}(M)$ of the multiplication group of $M$ and the commutator ideal
${\bf g}'$ of ${\bf g}$. Moreover, for all ideals ${\bf b}$ of ${\bf g}$ such that $\hbox{dim}({\bf b}) \ge 3$,
${\bf a} \subset {\bf b} \subseteq {\bf g}'$ the orbit $B(e)$, where $B$ is the simply connected Lie group of ${\bf b}$, coincides with $M$. One has
${\bf a} \cap {\bf inn(L)}=\{ 0 \}$ and the intersection ${\bf b} \cap {\bf inn(L)}$ has dimension
$\hbox{dim}({\bf b})-2$.

If the factor loop $L/A(e)$ is isomorphic to $\mathbb R^2$, then we have case (i).
\end{lemma}

\section{There does not exist any $3$-dimensional topological loop
having a $6$-dimensional solvable decomposable Lie group with discrete
centre as its multiplication groups}
We may assume that the loop $L$ is simply connected and hence it is homeomorphic to $\mathbb R^{3}$
(cf. Lemma \ref{simplyconnected}). As the multiplication group $Mult(L)$ of $L$ is a $6$-dimensional decomposable solvable Lie group with discrete centre
	for the Lie algebra ${\bf mult(L)}$ of $Mult(L)$ we have the following possibilities:
	${\bf l}_2 \oplus {\bf l}_2 \oplus {\bf l}_2$, ${\bf g}_{3,i} \oplus {\bf g}_{3,j}$,
	${\bf l}_2 \oplus {\bf g}_{4,k}$, where
	${\bf g}_{3,i}$, ${\bf g}_{3,j}$, $i,j \in \{2,3,4,5\}$, are the $3$-dimensional solvable Lie algebras with trivial centre (cf. § 4 in \cite{mubarakzyanov1},  p. 119), ${\bf g}_{4,k}$, $k=2,4,5,6,7,10$,
	${\bf g}_{4,8}^{h \neq -1}$, ${\bf g}_{4,9}^{p \neq 0}$ are the $4$-dimensional solvable Lie algebras with trivial centre (see § 5 in \cite{mubarakzyanov1}, pp. 120-121).
These Lie algebras have trivial centre and neither
a subalgebra nor a factor Lie algebra is isomorphic to an elementary filiform
Lie algebra ${\bf f}_n$, $n=4, 5$. 	

The Lie algebras
	${\bf mult(L)}= {\bf l}_2 \oplus {\bf g}_{4,k}$, $k=2,4,5,6,7,10$,
	${\bf l}_2 \oplus {\bf g}_{4,8}^{h \neq -1}$, ${\bf l}_2 \oplus {\bf g}_{4,9}^{p \neq 0}$, where
	${\bf l}_2=\langle f_1,f_2 \rangle $, have the $1$-dimensional ideal ${\bf i}=\langle f_1 \rangle $. There does not exist any ideal ${\bf s}$ of ${\bf mult(L)}$ such that
	${\bf i} \subseteq {\bf s}$ and ${\bf mult(L)}/{\bf s}$ is isomorphic to the Lie algebra
	${\bf l}_2 \oplus {\bf l}_2$. By Lemma \ref{lemmaregi} c) these Lie algebras are not the Lie algebra of the multiplication group of a loop $L$. 	
	
Now we consider the Lie algebras
${\bf g}_{i,j}={\bf g}_{3,i} \oplus {\bf g}_{3,j}=\langle e_1, e_2, e_3 \rangle \oplus \langle e_4, e_5, e_6 \rangle $,
$i,j \in \{2,3,4,5\}$.
Let be $j=5$.
The Lie algebra ${\bf g}_{3,5}$, respectively ${\bf g}_{4,5}$ is defined by $[e_1,e_3]=e_1$, $[e_2,e_3]=he_2$,
$[e_4,e_6]=pe_4-e_5$, $[e_5,e_6]=e_4+pe_5$, $p \ge 0$, where $h=1$, respectively $-1 \le h <1$,
whereas
the Lie algebra ${\bf g}_{2,5}$ is given by
$[e_1,e_3]=e_1$, $[e_2,e_3]=e_1+e_2$, $[e_4,e_6]=pe_4-e_5$, $[e_5,e_6]=e_4+pe_5$, $p \ge 0$. They have the $1$-dimensional ideal
	${\bf i}=\langle e_1 \rangle $. There does not exist any ideal ${\bf s}$ of ${\bf g}_{i,5}$, $i=2,3,4$, such that
	${\bf i} \subseteq {\bf s}$ and ${\bf g}_{i,5}/{\bf s}$ is isomorphic to the Lie algebra
	${\bf l}_2 \oplus {\bf l}_2$.
The Lie algebra
	${\bf g}_{5,5}$ defined by $[e_1,e_3]=p_1 e_1-e_2$, $[e_2,e_3]=e_1+p_1 e_2$, $[e_4,e_6]=p_2 e_4-e_5$,
	$[e_5,e_6]=e_4+p_2 e_5$ with $p_1,p_2\geq 0$ has the minimal ideals ${\bf s}_1=\langle e_1, e_2 \rangle $,
	${\bf s}_2=\langle e_4, e_5 \rangle $. Let $S_i$, $i=1,2$, be the simply connected Lie groups of ${\bf s}_i$. If
	${\bf g}_{5,5}$ would be the Lie algebra of the multiplication group of $L$, then by Lemma
	\ref{lemmaregi} b) and c) at least one of the orbits $S_i(e)$, $i=1,2$, is a normal subloop of $L$ isomorphic to $\mathbb R$.  For this orbit the factor loop $L/S_i(e)$ is isomorphic to the group ${\mathcal L}_2$.
Since the factor Lie algebras
	${\bf g}_{5,5}/{\bf s}_i$, $i=1,2$, are not isomorphic to the Lie algebra
	${\bf l}_2 \oplus {\bf l}_2$ the Lie algebra ${\bf g}_{5,5}$ is excluded (cf. Lemma \ref{tool} (ii)).

The Lie algebras ${\bf g}_{3,3}$, ${\bf g}_{3,4}$, ${\bf g}_{4,4}$ are defined by
$[e_1,e_3]=e_1$, $[e_2,e_3]=h_1 e_2$, $[e_4,e_6]=e_4$, $[e_5,e_6]=h_2 e_5$ such that for ${\bf g}_{3,3}$ one has $h_1=h_2=1$, for
${\bf g}_{3,4}$ we have $h_1=1$, $-1 \le h_2 <1$ and for ${\bf g}_{4,4}$ one has $-1 \le h_1,h_2 <1$.
The Lie algebra ${\bf g}_{2,3}$, respectively ${\bf g}_{2,4}$ is given by $[e_1,e_3]=e_1$,
$[e_2,e_3]=e_1+e_2$, $[e_4,e_6]=e_4$, $[e_5,e_6]=h_2 e_5$, where $h_2=1$, respectively $-1 \le h_2 <1$. The Lie algebra
${\bf g}_{2,2}$ is defined by $[e_1,e_3]=e_1$, $[e_2,e_3]=e_1+e_2$,
$[e_4,e_6]=e_4$, $[e_5,e_6]=e_4+e_5$.
All these Lie algebras have the ideals ${\bf i}_1=\langle e_1 \rangle $, ${\bf i}_2=\langle e_4 \rangle $. Additionally, the Lie algebra ${\bf g}_{3,3}$ has the ideals ${\bf i}_3=\langle e_2+ l_1 e_1 \rangle $,
${\bf i}_4=\langle e_5+ l_2 e_4 \rangle $, $l_1, l_2 \in \mathbb R$, the Lie algebra ${\bf g}_{4,4}$ has the ideals
${\bf i}_5=\langle e_2 \rangle $, ${\bf i}_6=\langle e_5 \rangle $, the Lie algebra ${\bf g}_{2,3}$ has the ideal ${\bf i}_4$,
the Lie algebra ${\bf g}_{2,4}$ has the ideal ${\bf i}_6$, and the Lie algebra ${\bf g}_{3,4}$ has the ideals
${\bf i}_3$, ${\bf i}_6$. All Lie algebras have the ideal
	${\bf s}_1=\langle e_1, e_4 \rangle $ containing ${\bf i}_1$, ${\bf i}_2$, such that the factor Lie algebras
	${\bf g}_{i,j}/{\bf s}_1$, $i,j \in \{2,3,4\}$ are isomorphic to ${\bf l}_2 \oplus {\bf l}_2$. Furthermore, the Lie algebra
	${\bf g}_{3,3}$ has
the ideal ${\bf s}_2=\langle e_2+ l_1 e_1, e_5+ l_2 e_4 \rangle $,  the Lie algebra ${\bf g}_{4,4}$ has the ideal
${\bf s}_3=\langle e_2, e_5 \rangle $, the Lie algebra ${\bf g}_{2,3}$ has the ideal
${\bf s}_4=\langle e_1, e_5+ l_2 e_4 \rangle $, 		
the Lie algebra ${\bf g}_{2,4}$ has the ideal ${\bf s}_5=\langle e_1, e_5 \rangle $ and the Lie algebra
${\bf g}_{3,4}$ has the ideal ${\bf s}_6=\langle e_2+ l_1 e_1, e_5 \rangle $ such that the factor Lie algebras
${\bf g}_{3,3}/{\bf s}_2$, ${\bf g}_{4,4}/{\bf s}_3$,
	 ${\bf g}_{2,3}/{\bf s}_4$, ${\bf g}_{2,4}/{\bf s}_5$, ${\bf g}_{3,4}/{\bf s}_6$ are isomorphic to
	${\bf l}_2 \oplus {\bf l}_2$.
	If ${\bf g}_{i,j}$, $i,j \in \{2,3,4\}$ is the Lie algebra of the multiplication group of a 3-dimensional
topological loop $L$, then the orbits $I_k(e)$, $k=1, \cdots ,6$, where $I_k=\exp({\bf i}_k)$ and $e$ is the identity
element of $L$, are $1$-dimensional normal subgroups of $L$ isomorphic to $\mathbb R$ and the factor loops $L/I_k(e)$ are isomorphic to ${\mathcal L}_2$ (cf. Lemma \ref{tool} (ii)).
All Lie algebras ${\bf g}_{i,j}$, $i,j \in \{2,3,4\}$, have the ideals ${\bf s}_7=\langle e_1, e_2 \rangle $,
${\bf s}_8=\langle e_4, e_5 \rangle $ such that the factor Lie algebras ${\bf g}_{i,j}/{\bf s}_l$, $l=7,8$, are not isomorphic to ${\bf l}_2 \oplus {\bf l}_2$. Hence the orbits $S_l(e)$,
where $S_l = \exp({\bf s}_l)$, $l=7,8$, and $e$ is the identity element of $L$, are $2$-dimensional normal subloops of $L$ and therefore one has ${\bf s}_l \cap {\bf inn}(L) = \{ 0\}$, $l=7,8$ (cf. Lemma \ref{tool}).
All Lie algebras ${\bf g}_{i,j}$, $i,j \in \{2,3,4\}$, have the commutator subalgebra
${\bf n}_1=\langle e_1, e_2, e_4, e_5 \rangle $. Their $5$-dimensional ideals are
${\bf v}_1=\langle e_1, e_2, e_4, e_5, e_3 \rangle $, ${\bf v}_{2,k}=\langle e_1, e_2, e_4, e_5, e_6+k e_3 \rangle $,
$k \in \mathbb R$. Denote by $N_1$ the simply connected Lie group of ${\bf n}_1$. By Lemma \ref{tool} b)  we have $N_1(e)=S_l(e)$, $l=7,8$. Therefore the intersection ${\bf n}_1 \cap  {\bf inn}(L)$ has dimension $2$. Hence the Lie algebra
${\bf inn}(L)$ has the basis elements $r_1=e_4+ a_1 e_1+ a_2 e_2$,
$r_2=e_5+b_1 e_1+b_2 e_2$ such that at least one of $a_1, a_2$ as well as $b_1, b_2$ are different from $0$ and
$a_1 b_2-a_2 b_1 \neq 0$.

All Lie algebras ${\bf g}_{i,j}$, $i,j \in \{2,3,4\}$, have the ideals ${\bf n}_2=\langle e_1, e_2, e_3 \rangle $,
${\bf n}_3=\langle e_4, e_5, e_6 \rangle $. As ${\bf s}_7 < {\bf n}_2$ and ${\bf s}_8 < {\bf n}_3$ the orbits
$N_j(e)$, where $N_j = \exp({\bf n}_j)$, $j=2, 3$, has dimension $2$ or $3$.
If $S_7(e)=N_2(e)$ or $S_8(e)=N_3(e)$, then one has $\hbox{dim}({\bf n}_2 \cap  {\bf inn}(L))=1$ or
$\hbox{dim}({\bf n}_3 \cap  {\bf inn}(L))=1$. Hence the Lie algebra ${\bf inn}(L)$ has the basis element either
$r_3=e_3+c_1 e_1+c_2 e_2$ or $r_3'=e_6+d_1 e_4+d_2 e_5$, $c_i, d_i \in \mathbb R$, $i=1,2$.
Since
$[r_1,r_3]$, respectively $[r_2,r_3']$ is a non-zero element of the ideal ${\bf s}_7$, respectively
${\bf s}_8$ the subspaces $\langle r_1, r_2, r_3 \rangle $, $\langle r_1, r_2, r_3' \rangle $ are not $3$-dimensional subalgebras of ${\bf g}_{i,j}$, $i,j \in \{2,3,4\}$. This contradiction gives that $N_2(e)=L$ and $N_3(e)=L$. As ${\bf n}_2 < {\bf v}_1$ and
${\bf n}_3 < {\bf v}_{2,0}$ we obtain that
$N_2(e)=V_1(e)=V_{2,0}(e)=N_3(e)=L$. By Lemma \ref{lemmaregi} c) there exists a parameter $k \in \mathbb R \setminus \{ 0 \}$ such that $V_{2,k}(e)$ is the $2$-dimensional normal subloop $S_7(e)=S_8(e)$. Hence one has
$\hbox{dim}({\bf v}_{2,k} \cap  {\bf inn}(L))=3$. Therefore the Lie algebra ${\bf inn}(L)$ has the basis element
$r_4=e_6+ k e_3+l_1 e_1+l_2 e_2$ for some $k \in \mathbb R \setminus \{ 0 \}$,
$l_i \in \mathbb R$, $i=1,2$.

The subspace $\langle r_1, r_2, r_4 \rangle $ is not a $3$-dimensional subalgebra of the Lie algebras
${\bf g}_{2,3}$, ${\bf g}_{2,4}$, ${\bf g}_{3,4}$. Hence these Lie algebras cannot be the Lie algebra of the group $Mult(L)$
of $L$.

The subspace $\langle r_1, r_2, r_4 \rangle $ forms a $3$-dimensional subalgebra of ${\bf g}_{2,2}$ if and only if
$k=1$, $a_2=0$ and $b_2=a_1 \neq 0$. Hence the subalgebra ${\bf inn}(L) < {\bf g}_{2,2}$ has the form
${\bf inn}(L)=\langle e_4+a_1 e_1, e_5+b_1 e_1+ a_1 e_2, e_6+e_3+l_1 e_1+ l_2 e_2 \rangle $, $a_1 \neq 0$,
$b_1, l_i \in \mathbb R$.  		

The subspace $\langle r_1, r_2, r_4 \rangle $ forms a $3$-dimensional subalgebra of ${\bf g}_{3,3}$ if and only if $k=1$. Hence 		
the subalgebra ${\bf inn}(L) < {\bf g}_{3,3}$ has the form
${\bf inn}(L)=\langle e_4+a_1 e_1+ a_2 e_2, e_5+b_1 e_1+ b_2 e_2, e_6+e_3+l_1 e_1+ l_2 e_2 \rangle $ such that at least one of
$a_1, a_2$ as well as $b_1, b_2$ are different from $0$ and
$a_1 b_2-a_2 b_1 \neq 0$.

The subspace $\langle r_1, r_2, r_4 \rangle $ forms a $3$-dimensional subalgebra of ${\bf g}_{4,4}$ if and only if one has either
$a_1=0=b_2$, $k=h_2=\frac{1}{h_1}$, or $a_2=0=b_1$, $k=1$, $h_2=h_1$. Therefore the subalgebra
${\bf inn}(L) < {\bf g}_{4,4}$ has either the form
${\bf inn}(L)=\langle e_4+a_2 e_2, e_5+b_1 e_1, e_6+k e_3+l_1 e_1+ l_2 e_2 \rangle $ such that $a_2 b_1 \neq 0$,
$k=h_2=\frac{1}{h_1}$, or ${\bf inn}(L)=\langle e_4+a_1 e_1, e_5+b_2 e_2, e_6+e_3+l_1 e_1+ l_2 e_2 \rangle $
such that $a_1 b_2 \neq 0$, $h_1=h_2$.

Using the automorphism
$\phi(e_1)=\frac{b_{2}e_{1}-a_{2}e_{2}}{a_{1}b_{2}-a_{2}b_{1}}$, $\phi(e_2)=\frac{b_{1}e_{1}-a_{1}e_{2}}{a_{2}b_{1}-a_{1}b_{2}}$,
$\phi(e_3)=e_3-l_1 \phi(e_1)-l_2 \phi(e_2)$, $\phi(e_i)=e_i$, $i=4,5,6$, of the Lie algebras ${\bf g}_{i,i}$, $i=2,3,4$, such that for ${\bf g}_{2,2}$ one has $a_2=0, b_2=a_1\neq 0$ and for ${\bf g}_{4,4}$
we have $a_2=b_1=0$, $h_2=h_1$, we can reduce
 ${\bf inn}(L)$ to ${\bf inn}(L)_{1}=\langle e_4+e_1, e_5+e_2, e_6+e_3\rangle $. Moreover, the automorphism
$\phi(e_1)=\frac{1}{b_1}e_1$, $\phi(e_2)=\frac{1}{a_2}e_2$,
$\phi(e_3)=e_3-\frac{h_{1}l_{1}}{b_1}e_1-\frac{h_{1}l_{2}}{a_2}e_2$, $\phi(e_i)=e_i$, $i=4,5,6$, of the Lie algebra
${\bf g}_{4,4}$ with $h_2=\frac{1}{h_1}$ reduces ${\bf inn}(L)$
  to ${\bf inn}(L)_{2}=\langle e_4+e_2, e_5+e_1, e_6+\frac{1}{h_1}e_3\rangle$. Linear representations of the simply connected Lie groups $G_{i,i}$,
$i=2,3,4$, are given as follows: for $G_{2,2}$ one has
        $$g(x_1,x_2,x_3,x_4,x_5,x_6) g(y_1,y_2,y_3,y_4,y_5,y_6)= $$
        $$ g(x_1+(y_1+x_{3}y_{2})e^{x_3},x_2+y_2e^{x_3},x_3+y_3,x_4+(y_4+x_{6}y_{5})e^{x_6},x_5+y_5e^{x_6},x_6+y_6),$$
for $G_{3,3}$, where $h_1=1$, and for $G_{4,4}$ with $h_2=h_1$ we have
        $$g(x_1+y_1e^{x_3},x_2+y_2e^{h_{1}x_{3}},x_3+y_3,x_4+y_4e^{x_6},x_5+y_5e^{h_{1}x_{6}},x_6+y_6),$$
for $G_{4,4}$, where $h_2=\frac{1}{h_1}$, one has
        $$ g(x_1+y_1e^{x_3},x_2+y_2e^{h_{1}x_3},x_3+y_3,x_4+y_4e^{x_6},x_5+y_5e^{\frac{x_6}{h_{1}}},x_6+y_6).$$
We get that the subgroup $Inn(L)_{1}$ of $G_{2,2}$, $G_{3,3}$ and $G_{4,4}$ with $h_2=h_1$ has the form
$Inn(L)_{1}=\{g(u_1,u_2,u_3,u_1,u_2,u_3); u_i \in \mathbb{R}\}$, $i=1,2,3$, and the subgroup
$Inn(L)_{2}$ of $G_{4,4}$ with $h_2=\frac{1}{h_1}$ is
$Inn(L)_{2}=\{g(u_2,u_1,\frac{1}{h_{1}}u_3,u_1,u_2,u_3); u_i \in \mathbb{R}\}$, $i=1,2,3$.
Two arbitrary left transversals to the groups $Inn(L)_{1}$ and $Inn(L)_{2}$ in $G_{i,i}$, $i=2,3,4$, are
$$A=\{g(u,v,w,f_1(u,v,w),f_2(u,v,w),f_3(u,v,w)): u,v,w \in \mathbb{R}\},$$
$$B=\{g(k,l,m,g_1(k,l,m),g_2(k,l,m),g_3(k,l,m)): k,l,m \in \mathbb{R}\},$$
where $f_{i}(u,v,w):\mathbb{R}^{3}\rightarrow \mathbb{R}$ and $g_{i}(k,l,m):\mathbb{R}^{3}\rightarrow \mathbb{R}$, $i=1,2,3$, are continuous functions with $f_{i}(0,0,0)=g_{i}(0,0,0)=0$.
For all $a\in A$, $b\in B$ the condition $a^{-1}b^{-1}ab\in Inn(L)_1$ holds if and only if
in the cases $G_{2,2}$, $G_{3,3}$ with $h_1=1$ and $G_{4,4}$ with $h_2=h_1$ the equation
\begin{equation} \label{g}le^{-h_{1}m}(1-e^{-h_{1}w})+ve^{-h_{1}w}(e^{-h_{1}m}-1) =\nonumber \end{equation}
\begin{equation} g_{2}(k,l,m)e^{-h_{1}g_3(k,l,m)}(1-e^{-h_{1}f_3(u,v,w)})+\nonumber \end{equation} 
\begin{equation} f_{2}(u,v,w)e^{-h_{1}f_3(u,v,w)}(e^{-h_{1}g_3(k,l,m)}-1), \end{equation}
and for $G_{2,2}$ the equation
\begin{equation} e^{-m}(1-e^{-w})(k-lm)+e^{-w}(e^{-m}-1)(u-vw)+(wl-mv)e^{-w-m}= \nonumber \end{equation}
\begin{equation}  e^{-g_3(k,l,m)}(1-e^{-f_3(u,v,w)})(g_{1}(k,l,m)-g_{2}(k,l,m)g_{3}(k,l,m))+ \nonumber \end{equation}
\begin{equation}  e^{-f_3(u,v,w)}(e^{-g_3(k,l,m)}-1)(f_1(u,v,w)-f_{2}(u,v,w)f_{3}(u,v,w))+ \nonumber \end{equation}
\begin{equation} \label{g22}(g_{2}(k,l,m)f_{3}(u,v,w)-f_{2}(u,v,w)g_{3}(k,l,m))e^{-f_3(u,v,w)-g_3(k,l,m)},\end{equation}
for $G_{3,3}$ with $h_1=1$ and for $G_{4,4}$ with $h_2=h_1$ the equation
\begin{equation} \label{g33}ke^{-m}(1-e^{-w})+ue^{-w}(e^{-m}-1)= \nonumber \end{equation}
\begin{equation}g_1(k,l,m)e^{-g_3(k,l,m)}(1-e^{-f_3(u,v,w)})+f_1(u,v,w)e^{-f_3(u,v,w)}(e^{-g_3(k,l,m)}-1) \end{equation}
are satisfied for all $k,l,m,u,v,w \in \mathbb R$.
The products $a^{-1}b^{-1}ab$ are contained in $Inn(L)_2$ if and only if the equations
\begin{equation}le^{-h_{1}m}(1-e^{-h_{1}w})+ve^{-h_{1}w}(e^{-h_{1}m}-1) = \nonumber \end{equation}
\begin{equation}\label{g44} g_1(k,l,m)e^{-g_3(k,l,m)}(1-e^{-f_3(u,v,w)})+f_1(u,v,w) e^{-f_3(u,v,w)}(e^{-g_3(k,l,m)}-1),
\end{equation}
\begin{equation}ke^{-m}(1-e^{-w})+ue^{-w}(e^{-m}-1)=\nonumber \end{equation}
\begin{equation} \label{g44second} g_{2}(k,l,m)e^{-\frac{1}{h_{1}}g_3(k,l,m)}(1-e^{-\frac{1}{h_{1}}f_3(u,v,w)})+\nonumber \end{equation}
\begin{equation}f_{2}(u,v,w)e^{-\frac{1}{h_{1}}f_3(u,v,w)}(e^{-\frac{1}{h_{1}}g_3(k,l,m)}-1) \end{equation}
are satisfied for all $u,v,w,k,l,m\in \mathbb{R}$. Equations (\ref{g}), respectively (\ref{g44}) is satisfied precisely if one has
$f_3(u,v,w)=w$, $f_2(u,v,w)=v$, $g_3(k,l,m)=m$, $g_2(k,l,m)=l$, respectively
$f_3(u,v,w)=h_{1}w$, $f_1(u,v,w)=v$, $g_3(k,l,m)=h_{1}m$, $g_1(k,l,m)=l$. Then $A\cup B$ does not generate the groups
$G_{i,i}$, $i=2,3$, $G_{4,4}$ with $h_2=h_1$ and
$G_{4,4}$ with $h_2=\frac{1}{h_1}$. By Proposition \ref{kepka} the Lie algebras ${\bf g}_{i,i}$, $i=2,3$, ${\bf g}_{4,4}$
with $h_2=h_1$ and with $h_2=\frac{1}{h_1}$, are not the Lie algebras of the groups $Mult(L)$ of $3$-dimensional topological
loops $L$.

Hence it remains to deal with the Lie algebra
${\bf g}={\bf l}_2 \oplus {\bf l}_2 \oplus {\bf l}_2=
\langle f_1,f_2 \rangle \oplus \langle f_3,f_4 \rangle \oplus \langle f_5,f_6 \rangle$ with the Lie brackets
$[f_1,f_2]=f_1$, $[f_3,f_4]=f_3$, $[f_5,f_6]=f_5$. The Lie algebra ${\bf g}$ has the $1$-dimensional ideals
${\bf i}_1=\langle f_1 \rangle $, ${\bf i}_2=\langle f_3 \rangle $, ${\bf i}_3=\langle f_5 \rangle$. The ideals
${\bf s}_1=\langle f_1,f_2 \rangle $, ${\bf s}_2=\langle f_3,f_4 \rangle $, ${\bf s}_3=\langle f_5,f_6 \rangle$
have the properties ${\bf i}_j \subset {\bf s}_j$ and ${\bf g}/{\bf s}_j$, $j=1,2,3$, are isomorphic
to ${\bf l}_2 \oplus {\bf l}_2$. If ${\bf g}$ is the Lie algebra of the multiplication group of $L$, then
the orbits $I_j(e)$, $j=1,2,3$, where $I_j$ is the simply connected Lie group of ${\bf i}_j$ and $e$ is the identity element
of $L$, are $1$-dimensional normal subloops of $L$ such that the factor loops $L/I_j(e)$ are isomorphic to the $2$-dimensional
non-abelian Lie group $\mathcal{L}_2$ (cf. Lemma \ref{tool} (ii)).

For the ideals ${\bf a}_1=\langle f_1,f_3 \rangle $, ${\bf a}_2=\langle f_1,f_5 \rangle $,
${\bf a}_3=\langle f_3,f_5 \rangle$ the factor Lie algebras ${\bf g}/{\bf a}_j$, $j=1,2,3$, are not isomorphic to
${\bf l}_2 \oplus {\bf l}_2$. Hence these ideals and the commutator ideal ${\bf g}'=\langle f_1, f_3, f_5 \rangle $ satisfy the condition of  Lemma \ref{tool} b). Therefore the orbits $A_j(e)$ and $G'(e)$, where
$A_j$, respectively $G'$ is the simply connected Lie group of ${\bf a}_j$, $j=1,2,3$, respectively ${\bf g}'$,
are the same $2$-dimensional normal subloop $M$ of $L$. Furthermore, one has
${\bf inn(L)}	\cap {\bf a}_j=\{0\}$ for all $j=1,2,3$ and $\hbox{dim}({\bf g}' \cap {\bf inn(L)})=1$.
The commutator subalgebra ${\bf inn(L)}'$ of ${\bf inn(L)}$ is the intersection ${\bf g}' \cap {\bf inn(L)}$.
  As every element of ${\bf inn(L)}'$ is contained in one of the ideals ${\bf a}_j$ and
${\bf inn(L)}	\cap {\bf a}_j=\{0\}$ for all $j=1,2,3$, the Lie algebra ${\bf inn(L)}$ is abelian.
The $5$-dimensional ideals of ${\bf g}$ are:
$${\bf v}_1=\langle f_1, f_3, f_5, f_2+k_1 f_6, f_4+k_2 f_6 \rangle , \ \ {\bf v}_2=\langle f_1, f_3, f_5,
f_2+k_3 f_4, f_6+k_4 f_4 \rangle , $$
$${\bf v}_3=\langle f_1, f_3, f_5, f_4+k_5 f_2, f_6+k_6 f_2 \rangle , \ \ k_i \in \mathbb R, \ i=1, \dots, 6. $$
 Each $3$-dimensional abelian subalgebra of a $5$-dimensional ideal ${\bf v}_j$, $j=1,2,3$, contains a non-trivial ideal of
${\bf g}$. Hence the Lie algebra ${\bf g}={\bf l}_2 \oplus {\bf l}_2 \oplus {\bf l}_2$ is not the Lie algebra of the multiplication group a $3$-dimensional topological loop.
This discussion yields the following:

\begin{proposition} \label{6dimdecomtrialcentre}
The $6$-dimensional decomposable solvable Lie algebras with trivial centre are not the Lie algebra of the multiplication group of a topological loop of dimension $3$.
\end{proposition}

\begin{corollary}
There does not exist any connected topological proper loop $L$ of dimension $\le 3$ having a solvable Lie group of dimension
$\le 6$ with discrete centre as the multiplication group of $L$.
\end{corollary}
\begin{proof}
A nilpotent multiplication Lie group has always non-discrete centre (cf. \cite{figula0}, Theorem 1 and \cite{figula2}, Theorem).
Let $\hbox{dim}(L)=3$. If the multiplication group of $L$ is solvable and has dimension $\le 5$, then it is decomposable
having $1$- or $2$-dimensional centre (see Propositions 12, 13, 14, 15, 17 in \cite{figula1}). For $6$-dimensional solvable Lie groups the assertion follows from Theorems 3.6, 3.7 in \cite{figula2}, Proposition 13 in \cite{figula3} and
Proposition \ref{6dimdecomtrialcentre}.
\end{proof}

\section{Topological loops having a
$6$-dimensional solvable decomposable Lie group with $1$-dimensional centre as their
multiplication group}
In this section we determine the $6$-dimensional decomposable solvable Lie groups with $1$-dimensional centre $Z$ which are the multiplication group $Mult(L)$ of a $3$-dimensional connected simply connected topological loop $L$. These loops have a centre
$Z(L) \cong \mathbb R$ such that the factor loop $L/Z(L)$ is isomorphic to $\mathbb R^2$.

\begin{proposition} \label{1dimcentre}
Let $L$ be a connected topological loop of dimension $3$ such that its multiplication group $Mult(L)$ is a $6$-dimensional decomposable solvable Lie group with $1$-dimensional centre. Then $L$ has nilpotency class $2$. Moreover, the following Lie algebra pairs can occur as the Lie algebra ${\bf g}$ of the group $Mult(L)$ and the subalgebra ${\bf k}$ of the subgroup $Inn(L)$:
\newline
\noindent
If ${\bf g}$ has the form
${\bf g}=\mathbb{R}\oplus {\bf h}=\langle f_1\rangle \oplus \langle e_1,e_2,e_3,e_4,e_5\rangle$, where ${\bf h}$
is a $5$-dimensional solvable indecomposable Lie algebra with trivial centre, then one has:
\begin{itemize}
  \item ${\bf g}_1=\mathbb{R}\oplus {\bf g}_{5,19}^{\alpha=0,\beta\neq 0}$: $[e_2,e_3]=e_1$, $[e_1,e_5]=e_1$, $[e_2,e_5]=e_2$, $[e_4,e_5]=\beta e_4$, ${\bf k}_{1,\epsilon}=\langle e_1+f_1, e_2+\epsilon f_1, e_4+f_1 \rangle$, $\epsilon=0,1$,
  \item ${\bf g}_2=\mathbb{R}\oplus {\bf g}_{5,20}^{\alpha=0}$: $[e_2,e_3]=e_1$, $[e_1,e_5]=e_1$, $[e_2,e_5]=e_2$,
	$[e_4,e_5]=e_1+e_4$, ${\bf k}_{2,\epsilon}=\langle e_1+f_1, e_2+\epsilon f_1, e_4+a_3 f_1\rangle$, $a_3\in \mathbb{R}$,
	$\epsilon=0,1$,
  \item ${\bf g}_3=\mathbb{R}\oplus {\bf g}_{5,27}$: $[e_2,e_3]=e_1$, $[e_1,e_5]=e_1$, $[e_3,e_5]=e_3+e_4$,
	$[e_4,e_5]=e_1+e_4$, ${\bf k}_{3}=\langle e_1+f_1, e_3, e_4+a_3f_1\rangle$, $a_3\in \mathbb{R}$,
  \item ${\bf g}_4=\mathbb{R}\oplus {\bf g}_{5,28}^{\alpha=0}$: $[e_2,e_3]=e_1$, $[e_1,e_5]=e_1$, $[e_3,e_5]=e_3+e_4$,
	$[e_4,e_5]=e_4$, ${\bf k}_{4}=\langle e_1+a_1f_1, e_3, e_4+f_1\rangle$, $a_1\in \mathbb{R}\backslash\{0\}$,
  \item ${\bf g}_5=\mathbb{R}\oplus {\bf g}_{5,32}$: $[e_2,e_4]=e_1$, $[e_3,e_4]=e_2$, $[e_1,e_5]=e_1$, $[e_2,e_5]=e_2$,
	$[e_3,e_5]=he_1+e_3$, ${\bf k}_{5}=\langle e_1+f_1, e_2+a_2 f_1, e_3\rangle$, $h, a_2\in \mathbb{R}$,
  \item ${\bf g}_6=\mathbb{R}\oplus {\bf g}_{5,33}$: $[e_1,e_4]=e_1$, $[e_3,e_4]=\beta e_3$, $[e_2,e_5]=e_2$,
	$[e_3,e_5]=\gamma e_3$, $\beta^{2}+\gamma^{2}\neq 0$, ${\bf k}_{6}=\langle e_1+f_1, e_2+f_1, e_3+f_1\rangle$,
  \item ${\bf g}_7=\mathbb{R}\oplus {\bf g}_{5,34}$: $[e_1,e_4]=\alpha e_1$, $[e_2,e_4]=e_2$, $[e_3,e_4]=e_3$,
	$[e_1,e_5]=e_1$, $[e_3,e_5]=e_2$, ${\bf k}_{7}=\langle e_1+f_1, e_2+f_1, e_3+a_3f_1\rangle$, $\alpha, a_3\in \mathbb{R}$,
  \item ${\bf g}_8=\mathbb{R}\oplus {\bf g}_{5,35}$: $[e_1,e_4]=he_1$, $[e_2,e_4]=e_2$, $[e_3,e_4]=e_3$,
	$[e_1,e_5]=\alpha e_1$, $[e_2,e_5]=-e_3$, $[e_3,e_5]=e_2$, $h^{2}+\alpha^{2}\neq 0$,
	${\bf k}_{8,1}=\langle e_1+f_1, e_2+f_1, e_3+a_3f_1\rangle$, $a_3\in \mathbb{R}$,
	${\bf k}_{8,2}=\langle e_1+f_1, e_2, e_3+f_1\rangle$.
\end{itemize}
\noindent
If ${\bf g}$ is the Lie algebra ${\bf l}_2 \oplus {\bf n}=\langle f_1,f_2 \rangle \oplus \langle e_1,e_2,e_3,e_4 \rangle$,
where ${\bf n}$ is a $4$-dimensional solvable Lie algebra with $1$-dimensional centre $\langle e_1 \rangle$, then we have:
\begin{itemize}
  \item ${\bf g}_9={\bf l}_{2}\oplus {\bf g}_{4,1}$: $[f_1,f_2]=f_1$, $[e_2,e_4]=e_1$, $[e_3,e_4]=e_2$,
	${\bf k}_{9}=\langle f_1+e_1, e_2+a_2e_1, e_3\rangle $, $a_2\in \mathbb{R}$,
  \item ${\bf g}_{10}={\bf l}_{2}\oplus {\bf g}_{4,3}$: $[f_1,f_2]=f_1$, $[e_1,e_4]=e_1$, $[e_3,e_4]=e_2$,
	${\bf k}_{10}=\langle f_1+e_2, e_1+e_2, e_3\rangle $.
\end{itemize}
\noindent
If ${\bf g}$ is one of the following Lie algebras ${\bf g}_{3,1}\oplus {\bf g}_{3,i}$ and
${\bf l}_2\oplus \mathbb{R}\oplus {\bf g}_{3,i}$, $i=2,3,4,5$, where ${\bf g}_{3,1}=\langle e_1,e_2,e_3\rangle$ is
the $3$-dimensional nilpotent Lie algebra having the centre $\langle e_1\rangle$ and
${\bf g}_{3,i}=\langle e_4,e_5,e_6\rangle$ is a $3$-dimensional solvable Lie algebra with trivial centre, then one has:
\begin{itemize}
  \item ${\bf g}_{11}={\bf g}_{3,1}\oplus {\bf g}_{3,2}$: $[e_2,e_3]=e_1$, $[e_4,e_6]=e_4$, $[e_5,e_6]=e_4+e_5$,
 ${\bf k}_{11,1}=\langle e_2, e_4+e_1, e_5\rangle$, ${\bf k}_{11,2}=\langle e_3, e_4+e_1, e_5\rangle$,
  \item ${\bf g}_{12}={\bf g}_{3,1}\oplus {\bf g}_{3,3}$: $[e_2,e_3]=e_1$, $[e_4,e_6]=e_4$, $[e_5,e_6]=e_5$,
	${\bf k}_{12,1}=\langle e_2, e_4+e_1, e_5+e_1\rangle$, ${\bf k}_{12,2}=\langle e_3, e_4+e_1, e_5+e_1 \rangle$,
  \item ${\bf g}_{13}={\bf g}_{3,1}\oplus {\bf g}_{3,4}$: $[e_2,e_3]=e_1$, $[e_4,e_6]=e_4$, $[e_5,e_6]=he_5$,
	$-1\leq h< 1$, $h\neq0$, ${\bf k}_{13,1}=\langle e_2, e_4+e_1, e_5+e_1\rangle$,
	${\bf k}_{13,2}=\langle e_3, e_4+e_1, e_5+e_1\rangle$,
  \item ${\bf g}_{14}={\bf g}_{3,1}\oplus {\bf g}_{3,5}$: $[e_2,e_3]=e_1$, $[e_4,e_6]=pe_4-e_5$, $[e_5,e_6]=e_4+pe_5$,
	$p\geq0$, ${\bf k}_{14,1}=\langle e_2, e_4+e_1, e_5+a_3e_1 \rangle$, ${\bf k}_{14,2}=\langle e_3, e_4+e_1, e_5+a_3e_1 \rangle$, $a_3\in \mathbb{R}\backslash\{0\}$, ${\bf k}_{14,3}=\langle e_2, e_4, e_5+e_1 \rangle$,
	${\bf k}_{14,4}=\langle e_3, e_4, e_5+e_1\rangle$,
  \item ${\bf g}_{15}={\bf l}_{2}\oplus \mathbb{R}\oplus {\bf g}_{3,2}$: $[f_1,f_2]=f_1$, $[e_4,e_6]=e_4$,
	$[e_5,e_6]=e_4+e_5$, ${\bf k}_{15}=\langle f_1+e_3, e_4+e_3, e_5\rangle$,
  \item ${\bf g}_{16}={\bf l}_{2}\oplus \mathbb{R}\oplus {\bf g}_{3,3}$: $[f_1,f_2]=f_1$, $[e_4,e_6]=e_4$,
	$[e_5,e_6]=e_5$, ${\bf k}_{16}=\langle f_1+e_3, e_4+e_3, e_5+e_3\rangle$,
  \item ${\bf g}_{17}={\bf l}_{2}\oplus \mathbb{R}\oplus {\bf g}_{3,4}$: $[f_1,f_2]=f_1$, $[e_4,e_6]=e_4$,
	$[e_5,e_6]=he_5$, $-1\leq h < 1$, $h\neq0$, ${\bf k}_{17}=\langle f_1+e_3, e_4+e_3, e_5+e_3\rangle$,
  \item ${\bf g}_{18}={\bf l}_{2}\oplus \mathbb{R}\oplus {\bf g}_{3,5}$: $[f_1,f_2]=f_1$, $[e_4,e_6]=pe_4-e_5$,
	$[e_5,e_6]=e_4+pe_5$, $p\geq0$, ${\bf k}_{18,1}=\langle f_1+e_3, e_4+e_3, e_5+a_3 e_3\rangle$, $a_3\in \mathbb{R}$,
	${\bf k}_{18,2}=\langle f_1+e_3, e_4, e_5+e_3\rangle$.
\end{itemize}
\end{proposition}

\begin{proof}
By Lemma \ref{simplyconnected} we may assume that the loop $L$ is simply connected and hence it is homeomorphic to
$\mathbb R^{3}$. Every $6$-dimensional decomposable solvable Lie algebra with
$1$-dimensional centre has one of the following forms: $\mathbb{R}\oplus {\bf h}$,
${\bf l}_{2}\oplus {\bf n}$, ${\bf g}_{3,1} \oplus {\bf g}_{3,i}$, and
${\bf l}_{2}\oplus \mathbb{R} \oplus {\bf g}_{3,i}$, where ${\bf h}$, ${\bf n}$, ${\bf g}_{3,i}$ are described in the assertion.  For ${\bf h}$ we have the following possibilities:
${\bf g}_{5,i} $, $i=7,9,11,12,13,16,17,18,21,23$,$ 24$, $27,31,32,33,34,35,36,37$ and
${\bf g}_{5,15}^{\gamma \neq 0} $, ${\bf g}_{5,j}^{\alpha=0}$, $j=19,20,28$, ${\bf g}_{5,k}^{p \neq 0} $,
$k=25,26$, ${\bf g}_{5,30}^{h \neq -2}$. For ${\bf n}$ one has the following Lie algebras ${\bf g}_{4,i}$,
$i=1,3$, ${\bf g}_{4,8}^{h=-1}$, ${\bf g}_{4,9}^{p=0}$ and for the $3$-dimensional solvable Lie algebras with trivial centre
we have ${\bf g}_{3,i}$, $i=2,3,4,5$ (cf.
\cite{mubarakzyanov1}, \S4, \S5, and \cite{mubarakzyanov2}, \S 10, p. 105-106).

To prove the first assertion we have to show that $L$ has a normal subloop $N$ isomorphic to $\mathbb R$ such that the factor loop $L/N$ is isomorphic to $\mathbb R^2$ (cf. Lemma \ref{lemmaregi} b) and d)).  Assume firstly that the Lie algebra of the multiplication group of $L$ has the form
$\mathbb{R}\oplus {\bf h}=\langle f_1 \rangle \oplus \langle e_1,e_2,e_3,e_4,e_5  \rangle$.
If ${\bf h}\neq {\bf g}_{5,i}$, $i=33,34$, then there does not exist any ideal ${\bf s}$ containing the centre
${\bf z}=\langle f_1 \rangle $ such that the factor Lie algebras $(\mathbb{R}\oplus {\bf h})/{\bf s}$ are isomorphic to
${\bf f}_n$, $n=4,5$ or to ${\bf l}_{2}\oplus {\bf l}_{2}$.
According to Lemma \ref{tool} a) the factor loop $L/Z(e)$, where
$Z=\exp({\bf z})$, is isomorphic to $\mathbb R^2$ and the orbit $Z(e)$ is the normal subloop $N$.

The Lie algebras $\mathbb{R}\oplus {\bf g}_{5,i}$, $i=33,34$, have no factor Lie algebras isomorphic to ${\bf f}_n$, $n=4,5$. The Lie algebra $\mathbb{R}\oplus {\bf g}_{5,34}$ has the $1$-dimensional ideal ${\bf i}=\langle e_1\rangle$. None of the factor Lie algebras $\mathbb{R}\oplus {\bf g}_{5,34}/{\bf s}$, where ${\bf s}$ is any ideal containing  ${\bf i}$, is isomorphic to
${\bf l}_{2}\oplus {\bf l}_{2}$. Therefore the orbit $I(e)$, where $I$ is the simply connected Lie group of ${\bf i}$, can choose as the normal subloop $N$.

The Lie algebra $\mathbb{R}\oplus {\bf g}_{5,33}$, $\beta^2+\gamma^2 \neq 0$, have the ideals
${\bf i}_1=\langle f_1 \rangle $, ${\bf i}_2=\langle e_1 \rangle $, ${\bf i}_3=\langle e_2 \rangle $,
${\bf i}_4=\langle e_3 \rangle $. If $\mathbb{R}\oplus {\bf g}_{5,33}$ is the Lie algebra of the multiplication group of $L$, then the orbits $I_j(e)$, $j \in \{1, 2, 3, 4\}$, are normal subgroups of $L$ isomorphic to $\mathbb R$. The factor loops $L/I_j(e)$,
$j \in \{1, 2, 3, 4\}$, are isomorphic either to
${\mathcal L}_2$ or to $\mathbb R^2$ (cf. Lemma \ref{tool} a).
If all factor loops $L/I_j(e)$, $j \in \{1, 2, 3, 4\}$, are isomorphic to ${\mathcal L}_2$, then by Lemma
\ref{tool} (ii) there are $2$-dimensional ideals ${\bf s}_j$,
$j \in \{1, 2, 3, 4\}$ such that ${\bf i}_j \subset {\bf s}_j$ and
the factor Lie algebras $\mathbb{R}\oplus {\bf g}_{5,33}/{\bf s}_j$ are isomorphic to
${\bf l}_2 \oplus {\bf l}_2$. For the ideal ${\bf s}_1={\bf s}_4=\langle f_1, e_3 \rangle $ one has
$\mathbb{R}\oplus {\bf g}_{5,33}/{\bf s}_l \cong {\bf l}_2 \oplus {\bf l}_2$, $l=1,4$. The factor Lie algebra
$\mathbb{R}\oplus {\bf g}_{5,33}/\langle f_1, e_1 \rangle $ is isomorphic to
${\bf l}_2 \oplus {\bf l}_2$ if and only if $\gamma=0$ and
$\mathbb{R}\oplus {\bf g}_{5,33}/\langle f_1, e_2 \rangle $ is isomorphic to
${\bf l}_2 \oplus {\bf l}_2$ precisely if $\beta=0$. This contradiction to
$\beta^2+\gamma^2 \neq 0$ yields that at least one of the factor loops $L/I_j(e)$, $j \in \{1, 2, 3, 4\}$, is isomorphic
to $\mathbb R^2$. For such $j \in \{1, 2, 3, 4\}$ the orbit $I_j(e)$ is the requested normal subgroup $N$ of $L$.

Hence $L$ is centrally nilpotent of class $2$.
By Lemma \ref{tool} (i) the Lie algebra $\mathbb{R}\oplus {\bf h}$ has a $4$-dimensional abelian ideal
${\bf p}={\bf z}\oplus {\bf k}$, where ${\bf z}=\langle f_1\rangle$ and ${\bf k}$ is the Lie algebra of the group $Inn(L)$ and
 ${\bf p}$ contains the commutator subalgebra of $\mathbb{R}\oplus {\bf h}$. According to Lemma
\ref{bruck} (3) the subalgebra ${\bf k}$ does not contain any non-zero ideal of ${\bf g}$ and the normalizer $N_{\bf g}({\bf k})$ of ${\bf k}$ in ${\bf g}$ is ${\bf p}$. The commutator subalgebra $\mathbb{R}\oplus {\bf h}$ coincides with the commutator subalgebra ${\bf h}'$ of ${\bf h}$. The intersection of ${\bf z}$ and ${\bf h}'$ is trivial. Since
${\bf h}'\subset {\bf p}$ the Lie algebra ${\bf h}$ has a $3$-dimensional abelian commutator subalgebra.
Then for the triples $({\bf g}, {\bf p}, {\bf k})$ we obtain:

(a) The Lie algebras $\mathbb{R}\oplus {\bf g}_{5,j}^{\alpha=0}$, $j=19,20$, have the ideal ${\bf p}=\langle f_1,e_1,e_2,e_4  \rangle$, the subalgebra ${\bf k}$ has the form: ${\bf k}_{a_1,a_2,a_3}=\langle e_1+a_1f_1, e_2+a_2f_1, e_4+a_3f_1\rangle$,
$a_i\in \mathbb{R}$, $i=1,2,3$ such that in the case $\mathbb{R}\oplus {\bf g}_{5,19}^{\alpha=0}$ one has $a_1\neq0, a_3\neq0$
since $\langle e_1\rangle$ and $\langle e_4\rangle$ are ideals of $\mathbb{R}\oplus {\bf g}_{5,19}^{\alpha=0}$,
\newline
in the case $\mathbb{R}\oplus {\bf g}_{5,20}^{\alpha=0}$ we have $a_1\neq0$ since $\langle e_1\rangle$ is an ideal of
$\mathbb{R}\oplus {\bf g}_{5,20}^{\alpha=0}$.
Applying the automorphism $\phi(f_{1})=f_{1}$, $\phi(e_{1})=a_1e_{1}$, $\phi(e_{2})=e_{2}$, $\phi(e_{3})=a_1e_{3}$,
$\phi(e_{4})=a_3e_{4}$, $\phi(e_{5})=e_{5}$ for the Lie algebra $\mathbb{R}\oplus {\bf g}_{5,19}^{\alpha=0}$ if $a_2=0$,  respectively $\phi(e_{2})=a_2e_{2}$, $\phi(e_{3})=\frac{a_1}{a_2}e_{3}$ if $a_2\neq0$ we can reduce ${\bf k}_{a_1,a_2,a_3}$ to
${\bf k}_{1,\epsilon}$, where $\epsilon$ is equal to $0$, respectively to $1$. Using the automorphism
$\phi(f_{1})=\frac{1}{a_1}f_{1}$, $\phi(e_{i})=e_{i}$, $i=1,2,3,4,5$ for the Lie algebra
$\mathbb{R}\oplus {\bf g}_{5,20}^{\alpha=0}$ if $a_2=0$,  respectively $\phi(f_{1})=f_{1}$, $\phi(e_{j})=a_1e_{j}$, $j=1,4$,
$\phi(e_{2})=a_2e_{2}$, $\phi(e_{3})=\frac{a_1}{a_2}e_{3}$, $\phi(e_{5})=e_{5}$ if $a_2\neq0$ the Lie algebra
${\bf k}_{a_1,a_2,a_3}$ reduces to ${\bf k}_{2,\epsilon}$, where $\epsilon$ is equal to $0$, respectively to $1$.

(b) For the Lie algebras $\mathbb{R}\oplus {\bf g}_{5,27}$ and $\mathbb{R}\oplus {\bf g}_{5,28}^{\alpha=0}$ we have ${\bf p}=\langle f_1,e_1,e_3,e_4  \rangle$, the subalgebra ${\bf k}$ has the form: ${\bf k}_{a_1,a_2,a_3}=\langle e_1+a_1f_1, e_3+a_2f_1, e_4+a_3f_1\rangle$, $a_i\in \mathbb{R}$, $i=1,2,3$ such that in the case $\mathbb{R}\oplus {\bf g}_{5,27}$ one has $a_1\neq0$, since $\langle e_1\rangle$ is an ideal of $\mathbb{R}\oplus {\bf g}_{5,27}$,
\newline
in the case $\mathbb{R}\oplus {\bf g}_{5,28}^{\alpha=0}$ one has $a_1 a_3\neq0$ since $\langle e_1\rangle$ and $\langle e_4\rangle$ are ideals of $\mathbb{R}\oplus {\bf g}_{5,28}^{\alpha=0}$. Using the automorphism $\phi(f_{1})=f_{1}$, $\phi(e_{i})=a_1e_{i}$,
$i=1,4$,  $\phi(e_{j})=e_{j}$, $j=2,5$, $\phi(e_{3})=a_1e_{3}+a_2e_1$ for $\mathbb{R}\oplus {\bf g}_{5,27}$, respectively
$\phi(e_{1})=a_1a_3e_{1}$, $\phi(e_{2})=a_1e_{2}$, $\phi(e_{3})=a_3e_{3}+a_2e_4$, $\phi(e_{4})=a_3e_{4}$ for
$\mathbb{R}\oplus {\bf g}_{5,28}^{\alpha=0}$ we can reduce ${\bf k}_{a_1,a_2,a_3}$ to ${\bf k}_{3}$, respectively
to ${\bf k}_{4}$ in the assertion.

(c) For the Lie algebras $\mathbb{R}\oplus {\bf g}_{5,i}$, $i=32,33,34,35$, have the ideal ${\bf p}=\langle f_1,e_1,e_2,e_3  \rangle$, and the subalgebra ${\bf k}$ has the form: ${\bf k}_{a_1,a_2,a_3}=\langle e_1+a_1f_1, e_2+a_2f_1, e_3+a_3f_1\rangle$,
$a_i\in \mathbb{R}$, $i=1,2,3$, such that in the case $\mathbb{R}\oplus {\bf g}_{5,32}$ we have $a_1\neq0$ since $\langle e_1\rangle$ is an ideal of $\mathbb{R}\oplus {\bf g}_{5,32}$,
\newline
in the case $\mathbb{R}\oplus {\bf g}_{5,33}$ we have $a_1a_2a_3\neq0$ since $\langle e_1\rangle$, $\langle e_2\rangle$ and
$\langle e_3\rangle$ are  ideals of $\mathbb{R}\oplus {\bf g}_{5,33}$,
\newline
in the case $\mathbb{R}\oplus {\bf g}_{5,34}$ we have $a_1a_2\neq0$ since $\langle e_1\rangle$, $\langle e_2\rangle$ are ideals of $\mathbb{R}\oplus {\bf g}_{5,34}$,
\newline
in the case $\mathbb{R}\oplus {\bf g}_{5,35}$ we have $a_1\neq0$ and at least one of $\{a_2,a_3\}$ is different from $0$ since
$\langle e_1\rangle$ and $\langle e_2, e_3\rangle$ are ideals of $\mathbb{R}\oplus {\bf g}_{5,35}$. The automorphism
$\phi(f_{1})=f_{1}$, $\phi(e_{i})=a_1e_{i}$, $i=1,2$, $\phi(e_{3})=a_1e_{3}+a_3e_1$ and $\phi(e_{j})=e_{j}$, $j=4,5$ for
$\mathbb{R}\oplus {\bf g}_{5,32}$, respectively $\phi(e_{2})=a_2e_{2}$, $\phi(e_{3})=a_3e_{3}$ for
$\mathbb{R}\oplus {\bf g}_{5,33}$, respectively $\phi(e_{s})=a_2e_{s}$, $s=2,3$ for $\mathbb{R}\oplus {\bf g}_{5,34}$ reduces the Lie algebra ${\bf k}_{a_1,a_2,a_3}$ to ${\bf k}_{5}$, respectively to ${\bf k}_{6}$, respectively to ${\bf k}_{7}$ in the assertion. Applying the automorphism $\phi(f_{1})=f_{1}$, $\phi(e_{1})=a_1e_{1}$, $\phi(e_{i})=a_2e_{i}$, $i=2,3$ and
$\phi(e_{j})=e_{j}$, $j=4,5$ for the Lie algebra $\mathbb{R}\oplus {\bf g}_{5,35}$ if $a_1 a_2\neq0$, respectively
$\phi(e_{s})=a_3e_{s}$, $s=2,3$ if $a_1 a_3\neq0$ and $a_2=0$ we can reduce ${\bf k}_{a_1,a_2,a_3}$ to ${\bf k}_{8,1}$, respectively ${\bf k}_{a_1,0,a_3}$ to ${\bf k}_{8,2}$ in the assertion.

Secondly, assume that the Lie algebra of the multiplication group of $L$ has the shape:
${\bf l}_{2}\oplus {\bf n}=\langle f_1,f_2 \rangle \oplus \langle e_1,e_2,e_3,e_4\rangle$ as in the assertion.
If ${\bf n}\neq {\bf g}_{4,1}$, then there does not exist any ideal ${\bf s}$ containing the ideal
${\bf i}=\langle f_1 \rangle $
such that the factor Lie algebras $({\bf l}_{2}\oplus {\bf n})/{\bf s}$ are isomorphic to ${\bf f}_n$, $n=4,5$ or to
${\bf l}_{2}\oplus {\bf l}_{2}$. The Lie algebra ${\bf l}_{2}\oplus {\bf g}_{4,1}$ has the centre ${\bf i}=\langle e_1\rangle$, but it has no factor Lie algebra isomorphic to ${\bf l}_{2}\oplus {\bf l}_{2}$. None of the factor Lie algebras
${\bf l}_{2}\oplus {\bf g}_{4,1}/{\bf s}$, where ${\bf s}$ is any ideal containing  ${\bf i}$, is isomorphic to ${\bf f}_n$,
$n=4,5$. Hence in both cases the orbit $I(e)$, where $I$ is the simply connected Lie group of ${\bf i}$, is a normal subgroup
of $L$ isomorphic to $\mathbb R$ with the property that the factor loop $L/I(e)$ is isomorphic to $\mathbb R^2$
(cf. Lemma \ref{tool} a). According to Lemma \ref{lemmaregi} a) and d) in both cases the orbit $I(e)$ coincides with the centre
$Z(L)$ of $L$ and $L$ has nilpotency class $2$.

Moreover, the Lie algebra ${\bf l}_{2}\oplus {\bf n}$ has a $4$-dimensional abelian ideal ${\bf p}={\bf z}\oplus {\bf k}$, where
${\bf z}$ is the $1$-dimensional centre of ${\bf l}_{2}\oplus \bf{n}$ and ${\bf k}$ is the Lie algebra of the group $Inn(L)$ such that ${\bf p}$ contains the commutator subalgebra of ${\bf l}_{2}\oplus {\bf n}$. The commutator subalgebra of
${\bf l}_{2}\oplus {\bf g}_{4,i}$, $i=1,3,8,9$, is the direct sum $\langle f_1 \rangle\oplus {\bf g}_{4,i}'$ where
${\bf g}_{4,i}'$ is the commutator subalgebra of ${\bf g}_{4,i}$. Since the commutator subalgebras ${\bf g}_{4,j}'$, $j=8,9$, are not abelian the Lie algebras ${\bf l}_{2}\oplus {\bf g}_{4,j}$, $j=8,9$, are excluded. Now we deal with the Lie algebras
${\bf l}_{2}\oplus {\bf g}_{4,k}$, $k=1,3$.

(d) The Lie algebra ${\bf l}_{2}\oplus {\bf g}_{4,1}$ has the centre ${\bf z}=\langle e_1 \rangle$ and ${\bf p}$ is the ideal
$\langle f_1,e_1,e_2,e_3  \rangle$. The subalgebra ${\bf k}$ has the form:
${\bf k}_{a_1,a_2,a_3}=\langle f_1+a_1e_1, e_2+a_2e_1, e_3+a_3e_1\rangle$, $a_i\in \mathbb{R}$, $i=1,2,3$, such that
$a_1\neq0$ since $\langle f_1\rangle$ is an ideal of
${\bf l}_{2}\oplus {\bf g}_{4,1}$. The centre of the Lie algebra ${\bf l}_{2}\oplus {\bf g}_{4,3}$ is
${\bf z}=\langle e_2 \rangle$ and the ideal ${\bf p}$ is again $\langle f_1,e_1,e_2,e_3  \rangle$. The subalgebra ${\bf k}$ has the form: ${\bf k}_{a_1,a_2,a_3}=\langle f_1+a_1e_2, e_1+a_2e_2, e_3+a_3e_2\rangle$, $a_i\in \mathbb{R}$, $i=1,2,3$ such that
$a_1\neq0$ and $a_2\neq 0$ since $\langle f_1\rangle$ and $\langle e_1\rangle$ are ideals of ${\bf l}_{2}\oplus {\bf g}_{4,3}$. The automorphism $\phi(f_{1})=a_1f_{1}$, $\phi(f_{2})=f_{2}$, $\phi(e_{3})=e_{3}-a_3e_1$, $\phi(e_{i})=e_{i}$, $i=1,2,4$ of
${\bf l}_{2}\oplus {\bf g}_{4,1}$, respectively $\phi(e_{1})=a_2e_{1}$, $\phi(e_{3})=e_{3}-a_3e_2$ of
${\bf l}_{2}\oplus {\bf g}_{4,3}$ reduces the Lie algebra
${\bf k}_{a_1,a_2,a_3}$ to ${\bf k}_{9}$, respectively to ${\bf k}_{10}$ in the assertion.

Finally, for the Lie algebras ${\bf g}_{3,1}\oplus {\bf g}_{3,i}=\langle e_1,e_2,e_3\rangle \oplus \langle e_4,e_5,e_6 \rangle$, respectively ${\bf l}_{2}\oplus \mathbb{R}\oplus {\bf g}_{3,i}=\langle f_1,f_2\rangle\oplus \langle e_3\rangle \oplus \langle e_4,e_5,e_6\rangle$, $i=2,3,4,5$, there does not exist any ideal ${\bf s}_1$, respectively ${\bf s}_2$ containing the ideal
${\bf i}_1=\langle e_1 \rangle $, respectively ${\bf i}_2=\langle f_1 \rangle $ such that the factor Lie algebras ${\bf g}_{3,1}\oplus {\bf g}_{3,i}/{\bf s}_1$, respectively, ${\bf l}_{2}\oplus \mathbb{R}\oplus {\bf g}_{3,i}/{\bf s}_2$ are isomorphic to
${\bf f}_n$, $n=4,5$ or to ${\bf l}_{2}\oplus {\bf l}_{2}$.
Hence if ${\bf g}_{3,1}\oplus {\bf g}_{3,i}$ or ${\bf l}_{2}\oplus \mathbb{R}\oplus {\bf g}_{3,i}$, $i=2,3,4,5$, is the Lie algebra of the multiplication group of $L$, then the orbits $I_i(e)$, $i=1,2$ are the centre of $L$ such that the factor loop
$L/I_i(e)$, $i=1,2$, are isomorphic to $\mathbb{R}^2$ (cf. Lemma \ref{lemmaregi} a), d). Hence $L$ is centrally nilpotent of class $2$.

According to Lemma \ref{tool} (i) we have to find an ideal ${\bf p}={\bf z}\oplus \bf{k} \cong \mathbb R^4$ of the Lie algebras
${\bf g}_{3,1}\oplus {\bf g}_{3,i}$ and ${\bf l}_{2}\oplus \mathbb{R}\oplus {\bf g}_{3,i}$, $i=2,3,4,5$, where ${\bf z}$ is
their $1$-dimensional centre, ${\bf p}$ contains their commutator subalgebra and ${\bf k}$ is the Lie algebra of the group
$Inn(L)$ satisfying the assertion of Lemma \ref{bruck} (3).

(e) The Lie algebras ${\bf g}_{3,1}\oplus {\bf g}_{3,i}$, $i=2,3,4,5$, have the centre ${\bf z}=\langle e_1 \rangle$ and the ideal ${\bf p}$ has one of the forms : ${\bf p}_{r}=\langle e_1,e_2+re_3,e_4,e_5 \rangle$, $r\in \mathbb{R}$, and ${\bf p}=\langle e_1,e_3,e_4,e_5 \rangle$. With respect to the ideals ${\bf p}_{r}$, ${\bf p}$ we obtain the subalgebras ${\bf k}_{r}=\langle e_2+re_3+a_1e_1, e_4+a_2e_1, e_5+a_3e_1 \rangle$, ${\bf k}_{a_1,a_2,a_3}=\langle e_3+a_1e_1 ,e_4+a_2e_1, e_5+a_3e_1 \rangle$, $r,a_i\in \mathbb{R}$, $i=1,2,3$, such that in the case ${\bf g}_{3,1}\oplus {\bf g}_{3,2}$ one has $a_2\neq0$ since $\langle e_4\rangle$ is an ideal of ${\bf g}_{3,1}\oplus {\bf g}_{3,2}$,
\newline
in the cases ${\bf g}_{3,1}\oplus {\bf g}_{3,i}$, $i=3,4$ we have $a_2a_3\neq0$  since $\langle e_4\rangle$, $\langle e_5\rangle$ are ideals of ${\bf g}_{3,1}\oplus {\bf g}_{3,i}$,
\newline
in the case ${\bf g}_{3,1}\oplus {\bf g}_{3,5}$ one has $a_2\neq0$ or $a_3\neq0$ since $\langle e_4, e_5\rangle$ is an ideal of ${\bf g}_{3,1}\oplus {\bf g}_{3,5}$. The automorphism $\phi(e_{2})=e_{2}-re_3-a_1e_1$, $\phi(e_{4})=a_2e_{4}$, $\phi(e_{5})=a_2e_{5}+a_3e_4$ and $\phi(e_{j})=e_{j}$, $j=1,3,6$, respectively $\phi(e_{2})=e_{2}$, $\phi(e_{3})=e_{3}-a_1e_1$ of ${\bf g}_{3,1}\oplus {\bf g}_{3,2}$ maps the subalgebra ${\bf k}_{r}$ onto ${\bf k}_{11,1}$, respectively ${\bf k}_{a_1,a_2,a_3}$ onto ${\bf k}_{11,2}$. The automorphism $\phi(e_{2})=e_{2}-re_3-a_1e_1$, $\phi(e_{4})=a_2e_{4}$, $\phi(e_{5})=a_3e_5$ and $\phi(e_{j})=e_{j}$, $j=1,3,6$, respectively $\phi(e_{2})=e_{2}$ and $\phi(e_{3})=e_{3}-a_1e_1$ of ${\bf g}_{3,1}\oplus {\bf g}_{3,i}$, $i=3,4$, maps the subalgebra ${\bf k}_{r}$ onto ${\bf k}_{12,1}={\bf k}_{13,1}$, respectively ${\bf k}_{a_1,a_2,a_3}$ onto ${\bf k}_{12,2}={\bf k}_{13,2}$ in the assertion.
For the Lie algebra ${\bf g}_{3,1}\oplus {\bf g}_{3,5}$ the automorphisms $\phi(e_{2})=e_2-re_3-a_1e_{1}$, $\phi(e_{j})=a_2e_{j}$, $j=4,5$, $\phi(e_{i})=e_{i}$, $i=1,3,6$, respectively $\phi(e_{2})=e_2$, $\phi(e_{3})=e_{3}-a_1e_1$ if $a_2\neq0$ reduce ${\bf k}_{r}$ to ${\bf k}_{14,1}$, respectively ${\bf k}_{a_1,a_2,a_3}$ to ${\bf k}_{14,2}$. Moreover, if $a_3\neq0$ and $a_2=0$, then the automorphisms $\phi(e_{2})=e_2-re_3-a_1e_{1}$, $\phi(e_{j})=a_3e_{j}$, $j=4,5$, $\phi(e_{i})=e_{i}$, $i=1,3,6$, respectively $\phi(e_{2})=e_2$, $\phi(e_{3})=e_{3}-a_1e_1$, change the Lie algebra ${\bf k}_{r}$ to ${\bf k}_{14,3}$, respectively ${\bf k}_{a_1,a_2,a_3}$ to ${\bf k}_{14,4}$ in the assertion.

The centre of the Lie algebras ${\bf l}_{2}\oplus \mathbb{R} \oplus {\bf g}_{3,i}$ with $i=2,3,4,5$, is
${\bf z}=\langle e_3\rangle$ and their ideal ${\bf p}$ is $\langle f_1,e_4,e_5,e_3\rangle$. The subalgebra ${\bf k}$ has the form: ${\bf k}_{a_1,a_2,a_3}=\langle f_1+a_1e_3, e_4+a_2e_3, e_5+a_3e_3\rangle$, $a_i\in \mathbb{R}$, $i=1,2,3$, such that in the case  ${\bf l}_{2}\oplus \mathbb{R} \oplus {\bf g}_{3,2}$ we have $a_1\neq0$ and $a_2\neq0$ since $\langle f_1\rangle$ and
$\langle e_4\rangle$ are  ideals of ${\bf l}_{2}\oplus \mathbb{R} \oplus {\bf g}_{3,2}$,
\newline
in the cases ${\bf l}_{2}\oplus \mathbb{R} \oplus {\bf g}_{3,i}$, $i=3,4$, one has $a_1a_2a_3\neq0$ since $\langle f_1\rangle$,
$\langle e_4\rangle$ and $\langle e_5\rangle$ are  ideals of ${\bf l}_{2}\oplus \mathbb{R} \oplus {\bf g}_{3,i}$,
\newline
in the case ${\bf l}_{2}\oplus \mathbb{R} \oplus {\bf g}_{3,5}$ we have $a_1\neq0$ and at least one of $\{a_2, a_3\}$ is different from $0$ since $\langle f_1\rangle$ and $\langle e_4, e_5\rangle$ are ideals of
${\bf l}_{2}\oplus \mathbb{R} \oplus {\bf g}_{3,5}$. Using the automorphism $\phi(f_{1})=a_1f_{1}$, $\phi(f_{2})=f_{2}$,
$\phi(e_{4})=a_2e_{4}$, $\phi(e_{5})=a_{2}e_{5}+a_3e_4$ and $\phi(e_{j})=e_{j}$, $j=3,6$, for
${\bf l}_{2}\oplus \mathbb{R} \oplus {\bf g}_{3,2}$, respectively $\phi(e_{5})=a_{3}e_{5}$ for
${\bf l}_{2}\oplus \mathbb{R} \oplus {\bf g}_{3,i}$, $i=3,4$, the Lie algebra ${\bf k}_{a_1,a_2,a_3}$ reduces to ${\bf k}_{15}$, respectively to ${\bf k}_{16}={\bf k}_{17}$. Applying the automorphism $\phi(f_{1})=a_1f_{1}$, $\phi(f_{2})=f_{2}$,
$\phi(e_{4})=a_2e_{4}$, $\phi(e_{5})=a_{2}e_{5}$ and $\phi(e_{j})=e_{j}$, $j=3,6$, for
${\bf l}_{2}\oplus \mathbb{R} \oplus {\bf g}_{3,5}$ if $a_1 a_2\neq0$, respectively $\phi(e_{4})=a_{3}e_{4}$ and
$\phi(e_{5})=a_{3}e_{5}$ if $a_1 a_3\neq0$ and $a_2=0$ we can reduce ${\bf k}_{a_1,a_2,a_3}$ to ${\bf k}_{18,1}$, respectively
${\bf k}_{a_1,0,a_3}$ to ${\bf k}_{18,2}$.
\end{proof}
Using  (\cite{strugar}, \S4) we obtain:

\begin{lemma}\label{Liegroups1}
The simply connected Lie group $G_i$ and its subgroup $K_i$ of the Lie algebra ${\bf g}_i$ and its subalgebra ${\bf k}_i$, $i=1,...,18$, is isomorphic to the linear group of matrices the multiplication of which is given by: for $i=1$
\begin{equation} g(x_1,x_2,x_3,x_4,x_5,x_6) g(y_1,y_2,y_3,y_4,y_5,y_6)= \nonumber \end{equation}
\begin{equation} g( x_1+(y_1-x_3y_2)e^{x_5},x_2+y_2e^{x_5},(x_3+y_3)e^{x_5+y_5},x_4+y_4e^{bx_5},x_5+y_5,x_6+y_6),\nonumber 
\end{equation}
\begin{equation} K_{1,\epsilon}=\{g(u_1,u_2,0,u_3,0,u_1+\epsilon u_2+u_3); u_i \in \mathbb R, i=1,2,3 \},  b\in \mathbb{R}\backslash\{0\}, \epsilon=0,1, \nonumber \end{equation}

for $i=2$
\begin{equation}g(x_1,x_2,x_3,x_4,x_5,x_6) g(y_1,y_2,y_3,y_4,y_5,y_6)=g(x_1+(y_1-x_3y_2+x_5y_4)e^{x_5}, \nonumber \end{equation}
\begin{equation} x_2+y_2e^{x_5},(x_3+y_3)e^{x_5+y_5},x_4+y_4e^{x_5},x_5+y_5,x_6+y_6), \nonumber \end{equation}
\begin{equation} K_{2,\epsilon}=\{g(u_1,u_2,0,u_3,0,u_1+\epsilon u_2+a_3u_3); u_i \in \mathbb R, i=1,2,3 \}, \epsilon=0,1, a_3\in \mathbb R,\nonumber \end{equation}

for $i=3$
\begin{equation} g(x_1,x_2,x_3,x_4,x_5,x_6) g(y_1,y_2,y_3,y_4,y_5,y_6)=  \nonumber \end{equation}
\begin{equation} g(x_1+(y_1+x_5y_4+\frac{1}{2}(2x_2+x_5^{2})y_3)e^{x_5},(x_2+y_2+x_5y_5+\frac{1}{2}y_5^{2}+\frac{1}{2}x_5^{2})e^{x_5+y_5},\nonumber \end{equation}
\begin{equation} x_3+y_3e^{x_5},x_4+(y_4+x_5y_3)e^{x_5},x_5+y_5,x_6+y_6), \nonumber \end{equation}
\begin{equation} K_{3}=\{g(u_1,0,u_2,u_3,0,u_1+a_3u_3); u_i \in \mathbb R, i=1,2,3 \}, a_3\in \mathbb R,\nonumber \end{equation}

for $i=4$
\begin{equation} g(x_1,x_2,x_3,x_4,x_5,x_6)g(y_1,y_2,y_3,y_4,y_5,y_6)= g(x_1+(y_1+x_2y_3)e^{x_5}, \nonumber \end{equation}
\begin{equation} (x_2+y_2)e^{x_5+y_5},x_3+y_3e^{x_5},x_4+(y_4+x_5y_3)e^{x_5},x_5+y_5,x_6+y_6), \nonumber \end{equation}
\begin{equation} K_{4}=\{g(u_1,0,u_2,u_3,0,a_1u_1+u_3); u_i \in \mathbb R, i=1,2,3 \}, a_1\in \mathbb{R}\backslash\{0\}, \nonumber \end{equation}

for $i=5$
\begin{equation} g(x_1,x_2,x_3,x_4,x_5,x_6) g(y_1,y_2,y_3,y_4,y_5,y_6)=\nonumber \end{equation}
\begin{equation} g(x_1+(y_1+x_4y_2+ax_5y_3+\frac{1}{2}x_4^{2}y_3)e^{x_5},x_2+(y_2+x_4y_3)e^{x_5}, \nonumber \end{equation}
\begin{equation} x_3+y_3e^{x_5},(x_4+y_4)e^{x_5+y_5},x_5+y_5,x_6+y_6), a\in \mathbb{R}, \nonumber \end{equation}
\begin{equation} K_{5}=\{g(u_1,u_2,u_3,0,0,u_1+a_2u_2);u_i \in \mathbb R, i=1,2,3 \}, a_2\in \mathbb R,\nonumber \end{equation}

for $i=6$
\begin{equation} g(x_1,x_2,x_3,x_4,x_5,x_6) g(y_1,y_2,y_3,y_4,y_5,y_6)= \nonumber \end{equation}
\begin{equation} g(x_1+y_1e^{x_4},x_2+y_2e^{x_5},x_3+y_3e^{ax_5+bx_4},x_4+y_4,x_5+y_5,x_6+y_6),  \nonumber \end{equation}
\begin{equation} K_{6}=\{g(u_1,u_2,u_3,0,0,u_1+u_2+u_3); u_i \in \mathbb R, i=1,2,3 \}, a^2+b^2\neq 0, \nonumber \end{equation}

for $i=7$
\begin{equation} g(x_1,x_2,x_3,x_4,x_5,x_6) g(y_1,y_2,y_3,y_4,y_5,y_6)=g(x_1+y_1e^{ax_4+x_5}, \nonumber \end{equation}
\begin{equation} x_2+(y_2+x_5y_3)e^{x_4},x_3+y_3e^{x_4},x_4+y_4e^{ax_4+x_5},(x_5+y_5)e^{x_4+y_4},x_6+y_6), \nonumber \end{equation}
\begin{equation} K_{7}=\{g(u_1,u_2,u_3,0,0,u_1+u_2+a_3u_3);u_i \in \mathbb R, i=1,2,3 \}, a, a_3\in \mathbb R,\nonumber \end{equation}

for $i=8$
\begin{equation} g(x_1,x_2,x_3,x_4,x_5,x_6) g(y_1,y_2,y_3,y_4,y_5,y_6)= g(x_1+y_1e^{ax_5+bx_4}, \nonumber \end{equation}
\begin{equation} x_2+(y_2cos (x_5)-y_3 sin (x_5))e^{x_4},x_3+(y_3cos (x_5)+y_2 sin (x_5))e^{x_4}, \nonumber \end{equation}
\begin{equation} x_4+y_4,x_5+y_5,x_6+y_6), a^2+b^2\neq 0,\nonumber \end{equation}
\begin{equation} K_{8,1}=\{g(u_1,u_2,u_3,0,0,u_1+u_2+a_3u_3);u_i \in \mathbb R, i=1,2,3 \}, a_3\in \mathbb R,\nonumber \end{equation}
\begin{equation} K_{8,2}=\{g(u_1,u_2,u_3,0,0,u_1+u_3); u_i \in \mathbb R, i=1,2,3 \}, \nonumber \end{equation}

for $i=9$
\begin{equation} g(x_1,x_2,x_3,x_4,x_5,x_6) g(y_1,y_2,y_3,y_4,y_5,y_6)= \nonumber \end{equation}
\begin{equation} g(x_1+y_1+x_4y_2+\frac{1}{2}x_4^{2}y_3,x_2+y_2+x_4y_3,x_3+y_3,x_4+y_4,x_5+y_5e^{x_6},x_6+y_6), \nonumber \end{equation}
\begin{equation} K_{9}=\{g(u_1+a_2u_2,u_2,u_3,0,u_1,0);u_i \in \mathbb R, i=1,2,3 \}, a_2\in \mathbb R,\nonumber \end{equation}

for $i=10$
\begin{equation} g(x_1,x_2,x_3,x_4,x_5,x_6) g(y_1,y_2,y_3,y_4,y_5,y_6)= \nonumber \end{equation}
\begin{equation} g(x_1+y_1e^{x_4},x_2+y_2+x_4y_3,x_3+y_3,x_4+y_4,x_5+y_5e^{x_6},x_6+y_6), \nonumber \end{equation}
\begin{equation} K_{10}=\{g(u_1,u_1+u_3,u_2,0,u_3,0); u_i \in \mathbb R, i=1,2,3 \}, \nonumber \end{equation}

for $i=11$
\begin{equation}g(x_1,x_2,x_3,x_4,x_5,x_6) g(y_1,y_2,y_3,y_4,y_5,y_6)=g(x_1+y_1+x_2y_3, \nonumber \end{equation}
\begin{equation} x_2+y_2,x_3+y_3,x_4+(y_4+x_6y_5)e^{x_6},x_5+y_5e^{x_6},x_6+y_6), \nonumber \end{equation}
\begin{equation} K_{11,1}=\{g(u_2,u_1,0,u_2,u_3,0); u_i \in \mathbb R, i=1,2,3 \}, \nonumber \end{equation}
\begin{equation} K_{11,2}=\{g(u_2,0,u_1,u_2,u_3,0); u_i \in \mathbb R, i=1,2,3 \}, \nonumber \end{equation}

for $i=12$
\begin{equation}g(x_1,x_2,x_3,x_4,x_5,x_6) g(y_1,y_2,y_3,y_4,y_5,y_6)= \nonumber \end{equation}
\begin{equation} g(x_1+y_1+x_2y_3,x_2+y_2,x_3+y_3,x_4+y_4e^{x_6},x_5+y_5e^{x_6},x_6+y_6), \nonumber \end{equation}
\begin{equation} K_{12,1}=\{g(u_2+u_3,u_1,0,u_2,u_3,0); u_i \in \mathbb R, i=1,2,3 \}, \nonumber \end{equation}
\begin{equation} K_{12,2}=\{g(u_2+u_3,0,u_1,u_2,u_3,0); u_i \in \mathbb R, i=1,2,3 \}, \nonumber \end{equation}

for $i=13$
\begin{equation}g(x_1,x_2,x_3,x_4,x_5,x_6) g(y_1,y_2,y_3,y_4,y_5,y_6)=g(x_1+y_1+x_2y_3, \nonumber \end{equation}
\begin{equation} x_2+y_2,x_3+y_3,x_4+y_4e^{x_6},x_5+y_5e^{hx_6},x_6+y_6), -1\leq h<1, h\neq 0, \nonumber \end{equation}
\begin{equation} K_{13,1}=\{g(u_2+u_3,u_1,0,u_2,u_3,0); u_i \in \mathbb R, i=1,2,3 \}, \nonumber \end{equation}
\begin{equation} K_{13,2}=\{g(u_2+u_3,0,u_1,u_2,u_3,0); u_i \in \mathbb R, i=1,2,3 \}, \nonumber \end{equation}

for $i=14$
\begin{equation}g(x_1,x_2,x_3,x_4,x_5,x_6) g(y_1,y_2,y_3,y_4,y_5,y_6)= g(x_1+y_1+x_2y_3,x_2+y_2,\nonumber \end{equation}
\begin{equation} x_3+y_3,x_4+(y_4cos (x_6)+y_5sin (x_6))e^{px_6},x_5+(y_5cos (x_6)-y_4sin (x_6))e^{px_6},x_6+y_6),  
\nonumber \end{equation}
\begin{equation} p\geq 0, K_{14,1}=\{g(u_2+a_3u_3,u_1,0,u_2,u_3,0); u_i \in \mathbb R, i=1,2,3, \}, a_3\in \mathbb{R}\backslash\{0\}, \nonumber \end{equation}
\begin{equation} K_{14,2}=\{g(u_2+a_3u_3,0,u_1,u_2,u_3,0); u_i \in \mathbb R, i=1,2,3, \}, a_3\in \mathbb{R}\backslash\{0\}, \nonumber \end{equation}
\begin{equation} K_{14,3}=\{g(u_3,u_1,0,u_2,u_3,0); u_i \in \mathbb R, i=1,2,3 \}, \nonumber \end{equation}
\begin{equation} K_{14,4}=\{g(u_3,0,u_1,u_2,u_3,0); u_i \in \mathbb R, i=1,2,3 \}, \nonumber \end{equation}

for $i=15$
\begin{equation} g(x_1,x_2,x_3,x_4,x_5,x_6) g(y_1,y_2,y_3,y_4,y_5,y_6)= \nonumber \end{equation}
\begin{equation} g(x_1+y_1e^{x_2},x_2+y_2,x_3+y_3e^{x_2},x_4+(y_4+x_6y_5)e^{x_6},x_5+y_5e^{x_6},x_6+y_6), \nonumber \end{equation}
\begin{equation} K_{15}=\{g(u_1,0,u_1+u_2,u_2,u_3,0); u_i \in \mathbb R, i=1,2,3 \}, \nonumber \end{equation}

for $i=16$
\begin{equation} g(x_1,x_2,x_3,x_4,x_5,x_6) g(y_1,y_2,y_3,y_4,y_5,y_6)= \nonumber \end{equation}
\begin{equation} g(x_1+y_1e^{x_2},x_2+y_2,x_3+y_3e^{x_2},x_4+y_4e^{x_6},x_5+y_5e^{x_6},x_6+y_6), \nonumber \end{equation}
\begin{equation} K_{16}=\{g(u_1,0,u_1+u_2+u_3,u_2,u_3,0); u_i \in \mathbb R, i=1,2,3 \}, \nonumber \end{equation}

for $i=17$
\begin{equation} g(x_1,x_2,x_3,x_4,x_5,x_6) g(y_1,y_2,y_3,y_4,y_5,y_6)=g(x_1+y_1e^{x_2}, \nonumber \end{equation}
\begin{equation} x_2+y_2,x_3+y_3e^{x_2},x_4+y_4e^{x_6},x_5+y_5e^{hx_6},x_6+y_6), -1\leq h <1, h\neq 0, \nonumber \end{equation}
\begin{equation} K_{17}=\{g(u_1,0,u_1+u_2+u_3,u_2,u_3,0); u_i \in \mathbb R, i=1,2,3 \}, \nonumber \end{equation}

for $i=18$
\begin{equation}g(x_1,x_2,x_3,x_4,x_5,x_6) g(y_1,y_2,y_3,y_4,y_5,y_6)= g(x_1+y_1e^{x_2},x_2+y_2,x_3+y_3e^{x_2}, \nonumber \end{equation}
\begin{equation}x_4+(y_4cos (x_6)+y_5sin (x_6))e^{px_6},x_5+(y_5cos (x_6)-y_4sin (x_6))e^{px_6},x_6+y_6), p\geq 0, \nonumber \end{equation}
\begin{equation}  K_{18,1}=\{g(u_1,0,u_1+u_2+a_3u_3,u_2,u_3,0);u_i \in \mathbb R, i=1,2,3 \}, a_3\in \mathbb R,\nonumber \end{equation}
\begin{equation} K_{18,2}=\{g(u_1,0,u_1+u_3,u_2,u_3,0); u_i \in \mathbb R, i=1,2,3 \}. \nonumber \end{equation}
\end{lemma}

\begin{proposition} \label{excluded1}
There does not exist any $3$-dimensional connected topological proper loop $L$ such that the Lie algebra ${\bf g}$ of the multiplication group of $L$ is one of the Lie algebras ${\bf g}_i$, $i=14,18$, with $p=0$.
\end{proposition}
\begin{proof}
We may assume that $L$ is simply connected and hence it is homeomorphic to $\mathbb R^3$ (cf. Lemma \ref{simplyconnected}). 
 We show that none of the groups
$G_i$, $i=14,18$, such that $p=0$ allows the existence of continuous left transversals $A$ and $B$ to $K_i$ in $G_i$ such that for all $a \in A$ and $b \in B$ one has $a^{-1} b^{-1} a b \in K_i$ and $A\cup B$ generates $G_i$. Hence Proposition
\ref{kepka} yields that the groups $G_i$, $i=14,18$, with $p=0$ are not the multiplication group of a loop
$L$. This proves the assertion.
\noindent
Two arbitrary left transversals to the groups $K_{14,i}$, $i=1,3$, in $G_{14}$ are:
\begin{equation} A =\{g(u, f_1(u,v,w),v,f_2(u,v,w) , f_3(u,v,w), w); u,v,w \in \mathbb R \},  \nonumber \end{equation}
\begin{equation} B =\{g(k, g_1(k,l,m), l, g_2(k,l,m), g_3(k,l,m), m); k,l,m \in \mathbb R \},  \nonumber \end{equation}
those to the groups $K_{14,j}$, $j=2,4$, in $G_{14}$ are:
\begin{equation} A =\{g(u,v,f_1(u,v,w),f_2(u,v,w) , f_3(u,v,w), w); u,v,w \in \mathbb R \},  \nonumber \end{equation}
\begin{equation} B =\{g(k,l,g_1(k,l,m), g_2(k,l,m), g_3(k,l,m), m); k,l,m \in \mathbb R \},  \nonumber \end{equation}
and those to the groups  $K_{18,j}$, $j=1,2$, in $G_{18}$ are:
\begin{equation} A =\{g(f_1(u,v,w),u,v,f_2(u,v,w), f_3(u,v,w), w); u,v,w \in \mathbb R \},  \nonumber \end{equation}
\begin{equation} B =\{g(g_1(k,l,m),k, l, g_2(k,l,m), g_3(k,l,m), m); k,l,m \in \mathbb R \},  \nonumber \end{equation}
where $f_i(u,v,w): \mathbb R^3 \to \mathbb R$ and $g_i(k,l,m): \mathbb R^3 \to \mathbb R$, $i =1,2,3$, are
continuous functions with $ f_i(0,0,0)= g_i(0,0,0)= 0$. The products $a^{-1}b^{-1}ab $, $a\in A$, $b \in B $ are elements
in $K_{14,1}$, respectively in $K_{14,2}$, respectively in $K_{18,1}$ if and only if
\newline
\begin{equation} (cos(m) -1)(f_2(u,v,w)(cos(w) +a_3 sin(w))+f_3(u,v,w)(a_3 cos(w)- sin(w)))+\nonumber \end{equation}
\begin{equation} (cos(w) -1)(g_3(k,l,m)(sin(m) -a_3 cos(m))-g_2(k,l,m)(cos(m) +a_3 sin(m)))- \nonumber \end{equation}
\begin{equation} sin(m)(f_2(u,v,w)(sin(w)-a_3cos(w))+f_3(u,v,w)(cos(w)+a_3sin(w)))+ \nonumber \end{equation}
\begin{equation} sin(w)(g_3(k,l,m)(cos(m)+a_3sin(m))+g_2(k,l,m)(sin(m)-a_3cos(m)))= \nonumber \end{equation}
\begin{equation} \label{g14}f_1(u,v,w)l-g_1(k,l,m)v, \end{equation}
respectively
\begin{equation} \label{g14one}g_1(k,l,m)v-f_1(u,v,w)l,\end{equation}
respectively
\begin{equation} \label{g18}e^{-u}(1-e^{-k})(f_1(u,v,w)-v)-e^{-k}(1-e^{-u})(g_1(k,l,m)-l) \end{equation}
are satisfied for all $k,l,m,u,v,w\in \mathbb{R}$, with $a_3 \in \mathbb{R}$. Moreover, the products $a^{-1}b^{-1}ab$,
$a\in A$, $b \in B $ are elements
in $K_{14,3}$, respectively in $K_{14,4}$, respectively in $K_{18,2}$ precisely if
\begin{equation} (cos(m)-1)(f_2(u,v,w)sin(w)+f_3(u,v,w)cos(w))-\nonumber \end{equation}
\begin{equation} (cos(w)-1)(g_2(k,l,m)sin(m)+g_3(k,l,m)cos(m))+\nonumber \end{equation}
\begin{equation} sin(m)(f_2(u,v,w)cos(w)-f_3(u,v,w)sin(w))+ \nonumber \end{equation} 
\begin{equation} sin(w)(g_3(k,l,m)sin(m)-g_2(k,l,m)cos(m))=\nonumber \end{equation}
\begin{equation} \label{g14second} f_1(u,v,w)l-g_1(k,l,m)v, \end{equation}
respectively
\begin{equation} \label{g14third}g_1(k,l,m)v-f_1(u,v,w)l,\end{equation}
respectively
\begin{equation} \label{g18second} e^{-u}(1-e^{-k})(f_1(u,v,w)-v)-e^{-k}(1-e^{-u})(g_1(k,l,m)-l) \end{equation}
hold for all $k,l,m,u,v,w \in \mathbb R$. The equations (\ref{g14}), (\ref{g14one}), (\ref{g18}), (\ref{g14second}), 
(\ref{g14third}) and (\ref{g18second}) are satisfied precisely if their left hand side as well as their right hand side are zero. The right hand side of these equations is zero if and only if $f_1(u,v,w)=v$ and $g_1(k,l,m)=l$. In that case the set $A\cup B$ does not generate $G_{14}$, respectively $G_{18}$.
\end{proof}

\begin{theorem} \label{multiplication1}
Let $L$ be a connected simply connected topological proper loop of dimension $3$ such that its multiplication group is a 
$6$-dimensional solvable decomposable Lie group having $1$-dimensional centre. Then the pairs of the Lie groups
$(G_i,K_i)$, $i=1,\cdots,18$, given in Lemma \ref{Liegroups1} such that for $i=14,18$ one has $p\neq 0$ are the multiplication group $Mult(L)$ and the inner mapping group $Inn(L)$ of $L$.
\end{theorem}
\begin{proof}
Taking into account Lemma \ref{1dimcentre} and Proposition \ref{excluded1} it remains to find for each group $G_i$,  
$i=1,\cdots,18$, in Lemma \ref{Liegroups1}, such that for $i=14,18$ one has $p\neq 0$, $K_i$-connected left transversals $A_i$, 
$B_i$ (cf. Proposition \ref{kepka}). 
The sets
\begin{equation} A_{1,0}=\{g(1-e^{v}-ue^{v}(1-e^{-bv}),e^{v}(1-e^{-bv}),u,ue^{bv-v},v,w); u, v, w \in \mathbb R\},  \nonumber \end{equation}
\begin{equation} B_{1,0}=\{g(1-e^{l}-ke^{l}(e^{-bl}-1),e^l(e^{-bl}-1),k,-ke^{bl-l},l,m); k, l, m \in \mathbb R\},\nonumber \end{equation}
respectively 
\begin{equation} A_{1,1}=\{g(1-e^{v}(2+u-e^{-bv}-ue^{-bv}),e^{v}(1-e^{-bv}),u,ue^{bv-v},v,w); u, v, w \in \mathbb R\}, \nonumber \end{equation}
\begin{equation} B_{1,1}=\{g(1-e^{l}(e^{-bl}+ke^{-bl}-k),e^l(e^{-bl}-1),k,-ke^{bl-l},l,m); k, l, m \in \mathbb R \} \nonumber \end{equation}
are $K_{1,0}$-, respectively $K_{1,1}$-connected left transversals in $G_{1}^{b\neq 0}$.
The set
\begin{equation} A_{2,0}=B_{2,0}=\{g(v^{2}-u^{2}-a_{3}v+e^{v}-1,u,u,v,v,w); u, v, w \in \mathbb R\}, \nonumber \end{equation}
respectively 
\begin{equation} A_{2,1}=B_{2,1}=\{g(v^{2}-u^{2}-u-a_{3}v+e^{v}-1,u,u,v,v,w); u, v, w \in \mathbb R\}, \nonumber \end{equation}
$a_3 \in \mathbb R $, is $K_{2,0}$, respectively $K_{2,1}$-connected left transversal in $G_2$. 
The set 
\begin{equation} A_{3}=B_{3}=\{g(e^{v}-1+(v-a_3)v(1+u+a_3v-\frac{1}{2}v^{2})+(u+a_3v-\frac{1}{2}v^{2})^{2},\nonumber \end{equation}
\begin{equation}u,u+a_3v-\frac{1}{2}v^{2},v(1+u+a_3v-\frac{1}{2}v^{2}),v,w); u,v,w \in \mathbb R\}, a_3 \in \mathbb R, \nonumber \end{equation}
respectively 
\begin{equation} A_4=B_4=\{g(-w,u,a_1u+v,1-e^{v}+a_1w+(a_1u+v)^{2},v,w); u, v, w \in \mathbb R\},  a_1\neq 0, \nonumber \end{equation} 
is $K_{3}$-connected, respectively $K_4$-connected left transversal in $G_3$, respectively in $G_4$.
The set $A_5=B_5=$ 
\begin{equation} \{g(1-e^{v}+2auv-\frac{1}{2}a_2u^{2}-aa_2v+ua_2^{2},av+\frac{1}{2}u^{2}-ua_2,u,u,v,w); u,v,w \in \mathbb R\},  \nonumber \end{equation}
$a_2\in \mathbb R$, is $K_5$-connected left transversal in $G_{5}^{a}$.
The sets
\begin{equation} A_6 =\{g(e^{u}-e^{u-av-bu},e^{v}-e^{v-u},e^{av+bu-v}-e^{av+bu},u,v,w); u,v,w \in \mathbb R\},\nonumber \end{equation}
\begin{equation} B_6 =\{g(e^{k}-e^{k-l},e^{l-al-bk}-e^{l},e^{al+bk}-e^{al+bk-k},k,l,m); k,l,m \in \mathbb R\} \nonumber \end{equation}
are $K_6$-connected left transversals in $G_{6}^{a^2+b^2\neq 0}$.
The set $A_7=B_7=$ 
\begin{equation} \{g(ve^{au+v-u},1-e^{u}-(a_3-v)(e^{u}-e^{u-au-v}),e^{u}-e^{u-au-v},u,v,w); u,v,w \in \mathbb R\}, \nonumber \end{equation}
$a_3 \in \mathbb R$, is $K_7$-connected left transversal in $G_{7}^{a}$.
The set 
\begin{equation} A_{8,1}=B_{8,1}=\{g(e^{av+bu-u}sin(v),  \nonumber \end{equation}
\begin{equation}\frac{1}{1+a_3^{2}}(e^{u}(1-e^{-av-bu})(sin(v)+a_3cos(v))+(e^{u}-cos(v))(cos(v)-a_3sin(v)),\nonumber  \end{equation}
\begin{equation}\frac{1}{1+a_{3}^{2}}((e^{u}-cos(v))(sin(v)+a_3cos(v))-(e^{u}-e^{u-av-bu})(cos(v)-a_3sin(v)), \nonumber  
\end{equation}
\begin{equation} u,v,w); u, v, w \in \mathbb R\}, a_3\in \mathbb R, \nonumber \end{equation}
respectively
\begin{equation} A_{8,2}=B_{8,2}=\{g(e^{av+bu-u}sin(v),(e^{u}-e^{u-av-bu})cos(v)-(e^{u}-cos(v))sin(v),\nonumber \end{equation}
\begin{equation} (e^{u}-e^{u-av-bu})sin(v)+(e^{u}-cos(v))cos(v),u,v,w); u, v, w \in \mathbb R\}, \nonumber \end{equation}
is $K_{8,1}$-, respectively $K_{8,2}$-connected left transversal in $G_{8}^{a^2+b^2\neq 0}$.
The sets
\begin{equation} A_{9}=\{g(u,v+(1-e^{-w})(a_2+v),1-e^{-w},v,-\frac{1}{2}v^{2}e^{w},w),u,v,w \in \mathbb R\},\nonumber \end{equation}
\begin{equation} B_{9}=\{g(k,l+(e^{-m}-1)(a_2+l),e^{-m}-1,l,\frac{1}{2}l^{2}e^{m},m),k,l,m \in \mathbb R\}, \nonumber \end{equation}
$a_2\in \mathbb R$, respectively 
\begin{equation} A_{10}=\{g(ve^{v},u,1-e^{-w},v,e^{w}-e^{w-v},w), u,v,w \in \mathbb R\}, \nonumber \end{equation}
\begin{equation} B_{10}=\{g(e^{l}-e^{l-m},k,e^{-l}-1,l,-le^{m},m), k,l,m \in \mathbb R\} \nonumber \end{equation}
are $K_9$-connected, respectively $K_{10}$-connected left transversals in $G_9$, respectively in $G_{10}$. 
The sets
\begin{equation} A_{11,1}=\{g(u,-we^{-w},v,e^{w}+vwe^{w}-1,ve^{w},w), u,v,w \in \mathbb R\}, \nonumber \end{equation}
\begin{equation} B_{11,1}=\{g(k,me^{-m},l,e^{m}-mle^{m}-1,-le^{m},m), k,l,m \in \mathbb R\}, \nonumber \end{equation}
respectively
\begin{equation} A_{11,2}=\{g(u,v,we^{-w},e^{w}+vwe^{w}-1,ve^{w},w), u,v,w \in \mathbb R\}, \nonumber \end{equation}
\begin{equation} B_{11,2}=\{g(k,l,-me^{-m},e^{m}-mle^{m}-1,-le^{m},m), k,l,m \in \mathbb R\}  \nonumber \end{equation}
are $K_{11,1}$-, respectively $K_{11,2}$-connected left transversals in $G_{11}$.
The sets
\begin{equation} A_{12,1}=\{g(u,e^{-w}-1,v,ve^{w}-u,u,w), u,v,w \in \mathbb R\}, \nonumber \end{equation}
\begin{equation} B_{12,1}=\{g(k,1-e^{-m},l,-le^{m}-k,k,m), k,l,m \in \mathbb R\}, \nonumber \end{equation}
respectively
\begin{equation} A_{12,2}=\{g(u,v,1-e^{-w},ve^{w}-u,u,w), u,v,w \in \mathbb R\}, \nonumber \end{equation}
\begin{equation} B_{12,2}=\{g(k,l,e^{-m}-1,-le^{m}-k,k,m), k,l,m \in \mathbb R\} \nonumber \end{equation}
are $K_{12,1}$-, respectively $K_{12,2}$-connected left transversals in $G_{12}$.
The sets
\begin{equation} A_{13,1}=\{g(u,1-e^{w},v,-ve^{w},e^{-w}-e^{-2w},w); u, v, w \in \mathbb R\},  \nonumber \end{equation}
\begin{equation} B_{13,1}=\{g(k,e^{-m}-1,l,e^{m}-e^{2m},le^{-m},m); k, l, m \in \mathbb R\}, \nonumber \end{equation}
respectively
\begin{equation} A_{13,2}=\{g(u,v,e^{w}-1,-ve^{w},e^{-w}-e^{-2w},w); u, v, w \in \mathbb R\},  \nonumber \end{equation}
\begin{equation} B_{13,2}=\{g(k,l,1-e^{-m},e^{m}-e^{2m},le^{-m},m); k, l, m \in \mathbb R\} \nonumber \end{equation}
are $K_{13,1}$-, respectively $K_{13,2}$-connected left transversals in $G_{13}^{h=-1}$ and the sets
\begin{equation} A_{13,3}=\{g(u,1-e^{-w},v,e^{w}-e^{w-hw},-ve^{hw},w); u, v, w \in \mathbb R\},  \nonumber \end{equation}
\begin{equation} B_{13,3}=\{g(k,e^{-hm}-1,l,le^{m},e^{hm}-e^{hm-m},m); k, l, m \in \mathbb R\}, \nonumber \end{equation}
respectively
\begin{equation} A_{13,4}=\{g(u,v,e^{-w}-1,e^{w}-e^{w-hw},-ve^{hw},w); u, v, w \in \mathbb R\},  \nonumber \end{equation}
\begin{equation} B_{13,4}=\{g(k,l,1-e^{-hm},le^{m},e^{hm}-e^{hm-m},m); k, l, m \in \mathbb R\} \nonumber \end{equation}
are $K_{13,1}$-, respectively $K_{13,2}$-connected left transversals in $G_{13}^{-1<h<1}$. 
The set 
\begin{equation} A_{14,1}=B_{14,1}=\{g(u,e^{-pw}sin(w), v,\nonumber \end{equation}
\begin{equation} \frac{1}{a_{3}^{2}+1}(e^{pw}v(sin(w)-a_3cos(w))+(cos(w)-e^{pw})(cos(w)+a_3sin(w))),\nonumber \end{equation}
\begin{equation}\frac{1}{a_{3}^{2}+1}(e^{pw}v(a_3sin(w)+cos(w))+(cos(w)-e^{pw})(a_3cos(w)-sin(w))),w);\nonumber \end{equation}
\begin{equation} u, v, w \in \mathbb R \}, \nonumber \end{equation}
respectively 
\begin{equation} A_{14,2}=B_{14,2}=\{g(u,v,-e^{-pw}sin(w), \nonumber \end{equation}
\begin{equation} \frac{1}{a_{3}^{2}+1}(e^{pw}v(sin(w)-a_3cos(w))+(cos(w)-e^{pw})(cos(w)+a_3sin(w))),\nonumber \end{equation}
\begin{equation}\frac{1}{a_{3}^{2}+1}(e^{pw}v(a_3sin(w)+cos(w))+(cos(w)-e^{pw})(a_3cos(w)-sin(w))),w); \nonumber \end{equation} 
\begin{equation} u, v, w \in \mathbb R\}, a_3\neq 0, \nonumber \end{equation} 
is $K_{14,1}$-, respectively $K_{14,2}$-connected left transversal in $G_{14}^{p\neq 0}$, and the set 
\begin{equation} A_{14,3}=B_{14,3}=\{g(u,e^{-pw}sin(w),v,sin(w)(cos(w)-e^{pw})-e^{pw}vcos(w),\nonumber \end{equation}
\begin{equation} e^{pw}vsin(w)+cos(w)(cos(w)-e^{pw}),w); u, v, w \in \mathbb R\},  \nonumber \end{equation}
respectively
\begin{equation} A_{14,4}=B_{14,4}=\{g(u,v,-e^{-pw}sin(w),sin(w)(cos(w)-e^{pw})-e^{pw}vcos(w),\nonumber \end{equation}
\begin{equation} e^{pw}vsin(w)+cos(w)(cos(w)-e^{pw}),w); u, v, w \in \mathbb R\}  \nonumber \end{equation}
is $K_{14,3}$-, respectively $K_{14,4}$-connected left transversal in $G_{14}^{p\neq 0}$.
The sets 
\begin{equation} A_{15}=\{g(e^{u-w}-e^{u}+v,u,v,e^{w}-e^{w-u}+w^{2},w,w), u,v,w \in \mathbb R\}, \nonumber \end{equation}
\begin{equation} B_{15}=\{g(e^{k}-e^{k-m}+l,k,l,e^{m-k}-e^{m}+m^{2},m,m), k,l,m \in \mathbb R\} \nonumber \end{equation}
are $K_{15}$-connected left transversals in $G_{15}$.
The set 
\begin{equation} A_{16}=B_{16}=\{g(e^{u}+v-1,u,v,e^{w}-u-1,u,w), u,v,w \in \mathbb R\} \nonumber \end{equation}
is $K_{16}$-connected left transversal in $G_{16}$.
The sets
\begin{equation} A_{17}=\{g(e^{u-hw}-e^{u}+v,u,v,e^{w}-e^{w-u},e^{hw}-e^{hw-w},w); u, v, w \in \mathbb R\},  \nonumber \end{equation}
\begin{equation} B_{17}=\{g(e^{k}-e^{k-m}+l,k,l,e^{m}-e^{m-hm},e^{hm-k}-e^{hm},m); k, l, m \in \mathbb R\} \nonumber \end{equation}
are $K_{17}$-connected left transversals in $G_{17}^{-1\leq h<1}$.
The set 
\begin{equation} A_{18,1}=B_{18,1}=\{g(e^{u}+v-1,u,v,\nonumber \end{equation}
\begin{equation} \frac{1}{1+a_3^{2}}((cos(w)-e^{pw})(a_3sin(w)+cos(w))-sin(w)(a_3cos(w)-sin(w)),\nonumber \end{equation}
\begin{equation}\frac{1}{1+a_{3}^{2}}((cos(w)-e^{pw})(a_3cos(w)-sin(w))+sin(w)(a_3sin(w)+cos(w)),w); \nonumber \end{equation}
\begin{equation} u, v, w \in \mathbb R\}, a_3\in \mathbb R, \nonumber \end{equation}
 respectively
\begin{equation} A_{18,2}=B_{18,2}=\{g(e^{u}+v-1,u,v,sin(w)(e^{pw}-cos(w))-cos(w)sin(w),\nonumber \end{equation}
\begin{equation} sin(w)^{2}+cos(w)(e^{pw}-cos(w)),w); u, v, w \in \mathbb R\},  \nonumber \end{equation}
is $K_{18,1}$-, respectively $K_{18,2}$-connected left transversal in $G_{18}^{p\neq0}$.

For all $i=1,\cdots, 18$, the set $A_i \cup B_i$ generates the group $G_i$. By Proposition \ref{kepka} the assertion is proved.  
\end{proof}

\section{Topological loops having a
$6$-dimensional solvable decomposable Lie group with $2$-dimensional centre as their 
multiplication group}

In this section we obtain the $6$-dimensional decomposable solvable Lie groups with $2$-dimensional centre $Z$ which can be represented as the multiplication group $Mult(L)$ of a $3$-dimensional connected simply connected topological proper loop $L$. 
These loops have a $2$-dimensional centre $Z(L)$ isomorphic to $\mathbb R^2$ such that the factor loop $L/Z(L)$ is isomorphic to 
$\mathbb R$.

\begin{proposition} \label{nilpotent}
Let $L$ be a connected simply connected topological proper loop of dimension $3$ such that its multiplication group is an at most
$6$-dimensional decomposable nilpotent Lie group. Then the loop $L$ is centrally nilpotent of class $2$ and the groups
$\mathbb R \times \mathcal{F}_4$, $\mathbb R \times \mathcal{F}_5$ are the multiplication group of $L$. 
\end{proposition}
\begin{proof}
Each nilpotent Lie group has a centre of dimension $\ge 1$. Hence, if the group $Mult(L)$ is decomposable and nilpotent, then it has a $2$-dimensional centre and the loop $L$ has nilpotency class $2$ (cf. Lemma \ref{lemmaregi}, a), e)). According to the list of the Lie algebras in \cite{mubarakzyanov1}, \S5, and
\cite{mubarakzyanov2}, p. 100, the Lie algebra of the group $Mult(L)$ is either the direct sum
${\bf g}_{3,1} \oplus {\bf g}_{3,1}$, where ${\bf g}_{3,1}$ is the $3$-dimensional nilpotent non-abelian Lie algebra, or 
$\mathbb R \oplus {\bf f}_n$, $n=4,5$, or
$\mathbb R \oplus {\bf g}_{5,i}$, $i=4,5,6$. By Lemma \ref{lemmaregi} e) the Lie algebra of $Mult(L)$ has a $5$-dimensional abelian ideal containing its centre and its commutator subalgebra. Since there does not exist any such ideal for the Lie algebras ${\bf g}_{3,1} \oplus {\bf g}_{3,1}$ and  $\mathbb R \oplus {\bf g}_{5,i}$, $i=4,5,6$, these Lie algebras are excluded. Now the assertion follows from Proposition 5.1. in \cite{figula6}, pp. 400-406.
\end{proof}

\begin{proposition}\label{2dimensionalcentre}
Let $L$ be a connected topological loop of dimension $3$ such that the Lie algebra ${\bf g}$ of its multiplication group is a 
$6$-dimensional decomposable solvable non-nilpotent Lie algebra with $2$-dimensional centre. Then $L$ is centrally nilpotent of class $2$. Moreover, the following Lie algebra pairs can occur as the Lie algebra ${\bf g}$ of the group $Mult(L)$ and the subalgebra $\bf{k}$ of the subgroup $Inn(L)$:
\newline
If ${\bf g}_i=\mathbb{R}^{2}\oplus {\bf n}_i=\langle f_1, f_2 \rangle \oplus \langle e_1, \cdots, e_4 \rangle$,
$i=1,\cdots,4$, where ${\bf n}$ is a
$4$-dimensional solvable indecomposable Lie algebra with trivial centre, then one has
\begin{itemize}
   \item ${\bf n}_{1}={\bf g}_{4,2}^{\alpha\neq 0}$: $[e_1,e_4]=\alpha e_1, [e_2,e_4]=e_2,
	[e_3,e_4]=e_2+e_3$, ${\bf k}_1=\langle e_{1}+f_{1},e_{2}+f_{1},e_{3}\rangle$,
   \item ${\bf n}_{2}={\bf g}_{4,4}$: $[e_1,e_4]=e_1, [e_2,e_4]=e_1+e_2,[e_3,e_4]=e_2+e_3$,
	${\bf k}_2=\langle e_{1}+f_{1},e_{2}+a_{2}f_{1},e_{3}+a_{3}f_{1}\rangle$, $a_2,a_3\in \mathbb{R}$, 
   \item ${\bf n}_{3}={\bf g}_{4,5}^{-1\leq\gamma\leq\beta\leq1,\gamma\beta\neq 0}$:
	$[e_1,e_4]=e_1, [e_2,e_4]=\beta e_2, [e_3,e_4]=\gamma e_3$, ${\bf k}_3=\langle e_{1}+f_{1},e_{2}+f_{1},e_{3}+f_{1}\rangle$, 
   \item ${\bf n}_{4}={\bf g}_{4,6}^{p\geq 0,\alpha\neq 0}$: $[e_1,e_4]=\alpha e_1, [e_2,e_4]=pe_2-e_3,
	[e_3,e_4]=e_2+pe_3$, ${\bf k}_{4,1}=\langle e_{1}+f_{1},e_{2}+f_{1},e_{3}+a_3f_1\rangle$, $a_3\in \mathbb{R}$,
	${\bf k}_{4,2}=\langle e_{1}+f_{1},e_{2},e_{3}+f_1\rangle$.
\end{itemize}
If 	${\bf g}_j=\mathbb{R}\oplus {\bf h}_j=\langle f_{1}\rangle\oplus\langle e_{1},e_{2}, e_{3},e_{4}, e_{5}\rangle$, where ${\bf h}_j$, $j=5, \cdots, 8$, is a $5$-dimensional solvable indecomposable Lie algebra with $1$-dimensional centre, then we have
	\begin{itemize}
   \item ${\bf h}_{5}={\bf g}_{5,8}^{0<|\gamma|\leq1}$: $[e_2,e_5]=e_1, [e_3,e_5]=e_3,[e_4,e_5]=\gamma e_4$, ${\bf k}_{5,\epsilon}=\langle e_{2}+\epsilon f_1,e_{3}+e_{1},e_{4}+e_{1}\rangle$, $\epsilon=0,1$, 
   \item ${\bf h}_{6}={\bf g}_{5,10}$: $[e_2,e_5]=e_1, [e_3,e_5]=e_2,[e_4,e_5]=e_4$, ${\bf k}_{6,\epsilon}=\langle e_{2},e_{3}+\epsilon f_1,e_{4}+e_{1}\rangle$, $\epsilon=0,1$, ${\bf k}_{6,3}=\langle e_{2}+b_1f_1,e_{3}+b_2f_1,e_{4}+f_{1}+ae_1\rangle$, $b_1, b_2 \in \mathbb R$, $a\neq 0$, 
   \item ${\bf h}_{7}={\bf g}_{5,14}^{p\neq0}$: $[e_2,e_5]=e_1, [e_3,e_5]=pe_3-e_4,[e_4,e_5]=e_3+pe_4$, ${\bf k}_{7,\epsilon}=\langle e_{2}+\epsilon f_1,e_{3}+e_{1},e_{4}+a_{3}e_{1}\rangle$, $\epsilon=0,1$, $a_3\in \mathbb{R}$, 
	${\bf k}_{7,\delta}=\langle e_{2}+\delta f_1,e_{3},e_{4}+e_{1}\rangle$, $\delta=0,1$, 
   \item ${\bf h}_{8}={\bf g}_{5,15}^{\gamma=0}$: $[e_1,e_5]=e_1, [e_2,e_5]=e_1+e_2, [e_4,e_5]=e_3$, ${\bf k}_{8,\epsilon}=\langle e_{1}+e_3,e_{2},e_{4}+\epsilon f_1\rangle$, $\epsilon=0,1$.
 \end{itemize}
\end{proposition}
\begin{proof}
By Lemma \ref{simplyconnected} we may assume that the loop $L$ is simply connected and hence it is homeomorphic to 
$\mathbb R^{3}$. As the multiplication group $Mult(L)$ of
$L$ is a $6$-dimensional decomposable solvable Lie group with $2$-dimensional centre the loop $L$ has nilpotency class
$2$ (cf. Lemma \ref{lemmaregi} a), e)). Furthermore, for the Lie algebra of $Mult(L)$ we have the following possibilities:
$\mathbb{R}^{2}\oplus {\bf n}$, $\mathbb{R}\oplus {\bf h} $, ${\bf l}_{2} \oplus \mathbb{R}\oplus {\bf g}_{3,1}$,  and
${\bf l}_{2}\oplus {\bf l}_{2}\oplus \mathbb{R}^{2}$, where ${\bf n}$ and ${\bf h}$ are characterized in the assertion. 
By \cite{mubarakzyanov1}, \S5, for ${\bf n}$ we have the Lie algebras ${\bf g}_{4,i}$, $i=2,4,5,6,7,10$,
${\bf g}_{4,8}^{h\neq -1}$, ${\bf g}_{4,9}^{p\neq 0}$. Moreover, the Lie algebras  ${\bf g}_{5,j}$,
$j=8,10,22,29,38,39$, ${\bf g}_{5,14}^{p\neq 0}$, ${\bf g}_{5,15}^{\gamma=0}$, 
${\bf g}_{5,19}^{\alpha=-1}$, ${\bf g}_{5,20}^{\alpha=-1}$, ${\bf g}_{5,28}^{\alpha=-1}$, 
${\bf g}_{5,25}^{p=0}$, ${\bf g}_{5,26}^{p=0}$ and ${\bf g}_{5,30}^{h=-2}$ can consider as ${\bf h}$ 
(cf. \cite{mubarakzyanov2}, \S 10, p. 105-106).

If these Lie algebras would be the Lie algebra of the multiplication group of $L$, then they have a $5$-dimensional abelian ideal containing their commutator ideal and their centre (cf. Lemma \ref{lemmaregi} e)). Since the Lie algebras  
${\bf l}_{2}\oplus \mathbb{R}\oplus {\bf g}_{3,1}$, ${\bf l}_{2}\oplus {\bf l}_{2} \oplus \mathbb{R}^{2}$, 
$\mathbb R^2 \oplus {\bf g}_{4,j}$, $j=7,10$, $\mathbb R^2 \oplus {\bf g}_{4,8}^{h \neq -1}$, 
$\mathbb R^2 \oplus {\bf g}_{4,9}^{p \neq 0}$, $\mathbb R \oplus {\bf g}_{5,r}^{\alpha=-1}$, $r=19,20,28$, 
$\mathbb R \oplus {\bf g}^{p=0}_{5,l}$, $l=25,26$, $\mathbb R \oplus {\bf g}_{5,30}^{h=-2}$, $\mathbb R \oplus {\bf g}_{5,p}$, 
$p=22,29,38,39$, do not contain any $5$-dimensional abelian ideal these Lie algebras are not the Lie algebra of the group 
$Mult(L)$ of $L$.
Hence it remain to deal with the Lie algebras ${\bf g}_{i}$, $i=1, \cdots, 8$ in the assertion.

The $1$-dimensional central subalgebras of ${\bf g}_{i}$, $i=1,2,3,4$, are
${\bf i}_{1}=\langle f_{2} \rangle$ and ${\bf i}_{2}=\langle f_{1}+a f_{2}\rangle$, $a\in \mathbb{R}$, those of ${\bf g}_{j}$, 
$j=5,6,7$, are
${\bf i}_{3}=\langle f_{1}+be_{1}\rangle$, $b \in \mathbb R$, and ${\bf i}_{4}=\langle e_{1}\rangle$, whereas those of
${\bf g}_{8}$ are ${\bf i}_{5}=\langle f_{1}+c e_{3}\rangle$, $c\in \mathbb{R}$, and
${\bf i}_{6}=\langle e_{3}\rangle$. With the exception of the Lie algebra ${\bf g}_{6}$ for every ideal ${\bf s}$ of each Lie algebra ${\bf g}_{i}$, $i=1,\cdots, 8$, such that
${\bf s}$ contains a $1$-dimensional central subalgebra of ${\bf g}_{i}$ the factor Lie algebras
${\bf g}_{i}/{\bf s}$ are not isomorphic to  ${\bf f}_{4}$. The Lie algebra ${\bf g}_{6}$ has the ideal
${\bf s}=\langle f_1+b e_1,e_4\rangle$ containing
${\bf i}_{3}$ such that the factor Lie algebra ${\bf g}_{6}/{\bf s}$ is isomorphic to ${\bf f}_{4}$.

According to Lemma \ref{lemmaregi} e) the simply connected Lie groups $G_i$ of ${\bf g}_{i}$, $i=1,\cdots,8$, has a $1$-dimensional connected central subgroup 
$N_d=\exp{{\bf n}_{d}}$, $d=1,2$, $4,\cdots,6$, such that the orbit $N_d(e)$ is isomorphic to $\mathbb{R}$ and the factor loop 
$L/N_d(e)$ is isomorphic to $\mathbb{R}^{2}$.
By Lemma \ref{tool} (i) the Lie algebras ${\bf g}_{i}$, $i=1, \cdots, 8$, have a $4$-dimensional abelian ideal
${\bf p}$ containing only a $1$-dimensional central subalgebra ${\bf n}_{d}$ of ${\bf g}_{i}$ and the Lie algebra ${\bf k}$ of 
the inner mapping group $Inn(L)$ of $L$ such that ${\bf g}_i' \subset {\bf p}$ and  ${\bf k}$ has the properties as in 
Lemma \ref{bruck} (3). Then for the triples $({\bf g}_i,{\bf p},{\bf k})$ we obtain:

(a) For the Lie algebras ${\bf g}_{i}$, $i=1,2,3,4$, the ideal ${\bf p}$ has one of the following forms 
${\bf p}_{a}=\langle f_{1}+af_{2},e_{1},e_{2},e_{3}\rangle$, $a \in \mathbb R$ and 
$\widetilde{{\bf p}}=\langle f_{2},e_{1},e_{2},e_{3}\rangle$. Hence for the subalgebras ${\bf k}$ one has  
${\bf k}_{a}=\langle e_{1}+a_{1}(f_{1}+af_{2}),e_{2}+a_{2}(f_{1}+af_{2}),e_{3}+a_{3}(f_{1}+af_{2})\rangle$, 
$a \in \mathbb R$ and $\widetilde{{\bf k}}=\langle e_{1}+a_{1}f_{2},e_{2}+a_{2}f_{2},e_{3}+a_{3}f_{2}\rangle$, where 
$a_i \in \mathbb{R}$, $i=1,2,3$. Using the automorphism $\phi(f_{1})=f_{2}$, $\phi(f_{2})=f_1+af_{2}$, 
$\phi(e_{i})=e_{i}$, $i=1,2,3,4$, the Lie algebra $\widetilde{{\bf k}}$ reduces to ${\bf k}_a$. So it remains to consider the subalgebra ${\bf k}_a$, such that
\newline
\noindent
in the case of the Lie algebra ${\bf g}_1$: $a_1a_2\neq0$ since $\langle e_{1}\rangle$ and $\langle e_{2}\rangle$ are ideals of 
${\bf g}_1$,
\newline
\noindent
in the case of ${\bf g}_2$: $a_1\neq0$ because $\langle e_{1}\rangle$ is an ideal of ${\bf g}_2$,
\newline
\noindent
in the case of ${\bf g}_3$: $a_1a_2a_3\neq0$ since $\langle e_{1}\rangle$, $\langle e_{2}\rangle$ and
$\langle e_{3}\rangle$ are ideals of ${\bf g}_3$,
\newline
\noindent
in the case of ${\bf g}_4$: $a_1\neq0$ and at least one of $\{a_2, a_3\}$ is different from $0$ because $\langle e_{1}\rangle$ and $\langle e_{2},e_{3}\rangle$ are ideals of ${\bf g}_4$.
\newline
\noindent
Using the automorphism $\phi(f_{1})=f_{1}-af_{2}$, $\phi(f_{2})=f_{2}$, $\phi(e_{1})=a_{1}e_{1}$, $\phi(e_{2})=a_{2}e_{2}$, $\phi(e_{3})=a_{2}e_{3}+a_{3}e_{2}$ and $\phi(e_{4})=e_{4}$ for ${\bf g}_1$, respectively $\phi(e_{j})=a_{1}e_{j}$, $j=2,3$ for ${\bf g}_2$, respectively $\phi(e_{3})=a_{3}e_{3}$ for ${\bf g}_3$, the Lie algebra ${\bf k}_a$ reduces to ${\bf k}_1$, respectively ${\bf k}_2$, $a_2,a_3\in \mathbb R$, respectively ${\bf k}_3$ in the assertion. Applying the automorphism $\phi(f_{1})=f_{1}-af_{2}$, $\phi(f_{2})=f_{2}$, $\phi(e_{1})=a_{1}e_{1}$, $\phi(e_{j})=a_{2}e_{j}$, $j=2,3$ and $\phi(e_{4})=e_{4}$ for ${\bf g}_4$ if $a_2\neq0$, respectively $\phi(e_{j})=a_{3}e_{j}$, $j=2,3$ if $a_2=0$ and $a_3\neq0$ we can reduce
${\bf k}_a$ to ${\bf k}_{4,1}$, $a_3\in \mathbb{R}$, respectively to ${\bf k}_{4,2}$ in the assertion.

(b) For the Lie algebras ${\bf g}_j$, $j=5,7$, the ideal $\bf{p}$ has one of the following shapes 
${\bf p}_a=\langle e_{1},f_1+a e_{2},e_{3},e_{4}\rangle$, $a\in \mathbb{R}\backslash \{0\}$, 
$\widetilde{{\bf p}}=\langle e_{1},e_{2},e_{3},e_{4}\rangle$. Hence the subalgebras
${\bf k}$ are ${\bf k}_a=\langle f_1+a e_{2}+a_{1}e_{1},e_{3}+a_{2}e_{1},e_{4}+a_{3}e_{1}\rangle$, $a\in \mathbb{R}\backslash\{0\}$, and $\widetilde{{\bf k}}=\langle e_{2}+a_{1}e_{1},e_{3}+a_{2}e_{1},e_{4}+a_{3}e_{1}\rangle$,
 $a_i\in \mathbb{R}$, $i=1,2,3$, such that
\newline
\noindent
for ${\bf g}_5$: $a_2a_3\neq 0$ since $\langle e_{3}\rangle$, $\langle e_{4}\rangle$ are ideals of ${\bf g}_5$, and
\newline
\noindent
for ${\bf g}_7$: $a_2\neq0$ or $a_3\neq0$ because $\langle e_{3}, e_{4}\rangle$ is an ideal of ${\bf g}_7$.
\newline
\noindent
The automorphism $\phi(f_{1})=f_{1}$, $\phi(e_{i})=e_{i}$, $i=1,5$, $\phi(e_{2})=e_{2}-a_{1}e_{1}$,
$\phi(e_{3})=a_{2}e_{3}$ and $\phi(e_{4})=a_{3}e_{4}$, respectively $\phi(f_{1})=af_{1}-a_1e_1$, $\phi(e_{2})=e_{2}$ of ${\bf g}_5$ maps the subalgebra $\widetilde{{\bf k}}$ onto ${\bf k}_{5,0}$, respectively the subalgebra
${\bf k}_a$ onto ${\bf k}_{5,1}$ in the assertion.
If $a_2\neq0$, then the automorphism $\phi(f_{1})=f_{1}$, $\phi(e_{i})=e_{i}$, $i=1,5$, $\phi(e_{2})=e_{2}-a_{1}e_{1}$, $\phi(e_{j})=a_{2}e_{j}$, $j=3,4$, respectively $\phi(f_{1})=af_{1}-a_1e_1$, $\phi(e_{2})=e_{2}$ of the Lie algebra
${\bf g}_7$ reduces the subalgebra $\widetilde{{\bf k}}$ to ${\bf k}_{7,\epsilon}$ with $\epsilon=0$, respectively the subalgebra
${\bf k}_a$ to ${\bf k}_{7,\epsilon}$ with $\epsilon=1$ in the assertion.
If $a_2=0$ and $a_3\neq0$, then using the automorphism $\phi(f_{1})=f_{1}$, $\phi(e_{i})=e_{i}$, $i=1,5$, 
$\phi(e_{2})=e_{2}-a_{1}e_{1}$, $\phi(e_{j})=a_{3}e_{j}$, $j=3,4$, respectively $\phi(f_{1})=af_{1}-a_1e_1$, $\phi(e_{2})=e_{2}$ of
${\bf g}_7$ we can change the subalgebra $\widetilde{{\bf k}}$ to ${\bf k}_{7,\delta}$ with $\delta=0$, respectively the subalgebra ${\bf k}_a$ to ${\bf k}_{7,\delta}$ such that $\delta=1$ in the assertion.

(c) For the Lie algebra ${\bf g}_8$ the ideal ${\bf p}$ has one of the following forms 
$\widetilde{{\bf p}}=\langle e_{1},e_{2},e_{3},e_{4}\rangle$, 
${\bf p}_a=\langle e_{1},e_{2},e_{3},f_1+a e_{4}\rangle$, $a\in \mathbb{R}\backslash \{0\}$. Therefore for the subalgebras
${\bf k}$ one has $\widetilde{{\bf k}}=\langle e_{1}+a_{1}e_{3},e_{2}+a_{2}e_{3},e_{4}+a_{3}e_{3}\rangle$,
${\bf k}_a=\langle e_{1}+a_{1}e_{3},e_{2}+a_{2}e_{3},f_{1}+a e_4+a_{3}e_{3}\rangle$,
$a\in \mathbb{R}\backslash\{0\}, a_i\in \mathbb{R}$, $i=1,2,3$, such that $a_1\neq 0$ since $\langle e_{1}\rangle$ is an ideal 
of ${\bf g}_8$. The automorphism $\phi(f_{1})=f_{1}$, $\phi(e_{i})=e_{i}$, $i=3,5$, $\phi(e_{1})=a_1e_1$,
$\phi(e_{2})=a_1e_{2}-a_{2}e_{1}$, and $\phi(e_{4})=e_4-a_{3}e_{3}$, respectively $\phi(f_{1})=af_{1}-a_3e_3$,
$\phi(e_{4})=e_4$, maps the subalgebra
$\widetilde{{\bf k}}$ onto ${\bf k}_{8,0}$, respectively ${\bf k}_a$ onto ${\bf k}_{8,1}$ of the assertion.

(d) If the Lie algebra ${\bf g}_6$ is the Lie algebra of the group $Mult(L)$ of $L$, then the factor loop $L/I_4(e)$, where 
$I_4=\exp ({\bf i}_4)$, is isomorphic to
$\mathbb R^2$. Hence the Lie algebra ${\bf k}$ of the group $Inn(L)$ of $L$ is a subalgebra of the ideal ${\bf p}$ having one of the following forms $\widetilde{{\bf p}}=\langle e_{1},e_{2},e_{4},e_{3}\rangle$,
${\bf p}_a=\langle e_{1},e_{2},e_{4},f_1+ae_{3}\rangle$, $a\in \mathbb R\backslash \{0\}$. Therefore we obtain the subalgebras
$\widetilde{{\bf k}}=\langle e_{2}+a_{1}e_{1},e_{3}+a_{2}e_{1},e_{4}+a_{3}e_{1}\rangle$,
${\bf k}_a=\langle e_{2}+a_{1}e_{1},f_1+a e_{3}+a_{2}e_{1},e_{4}+a_{3}e_{1}\rangle$, where
$a\in \mathbb{R}\backslash\{0\}, a_i\in \mathbb R$, $i=1,2$, and $a_3\neq 0$ since $\langle e_{4}\rangle$ is an ideal of
${\bf g}_6$. With the automorphism
$\phi(f_{1})=f_{1}$, $\phi(e_{i})=e_{i}$, $i=1,5$, $\phi(e_{2})=e_{2}-a_{1}e_{1}$, $\phi(e_{3})=e_3-a_{1}e_{2}-a_2e_1$ and
$\phi(e_{4})=a_{3}e_{4}$, respectively $\phi(f_{1})=af_{1}-a_2e_1$, $\phi(e_{3})=e_3-a_{1}e_{2}$, we can change the subalgebra 
$\widetilde{{\bf k}}$ onto ${\bf k}_{6,0}$, respectively
${\bf k}_a$ onto ${\bf k}_{6,1}$ in the assertion.

Since for the ideal
${\bf s}=\langle f_{1}+be_{1}, e_4\rangle$, $b\in \mathbb{R}$, of ${\bf g}_6$, the factor Lie algebra
${\bf g}_6/{\bf s}$ is isomorphic to ${\bf f}_4$ the factor loop $L/I_{3}(e)$, where 
$I_3=\exp ({\bf i}_3)$, is isomorphic to an elementary filiform loop $L_{\mathcal F}$. The orbit $S(e)$, where $S=exp(s)$, coincides with $I_{3}(e)$ (cf. Lemma \ref{lemmaregi} e). Hence the Lie algebra
${\bf k}$ contains the basis element $e_4+a_3(f_1+ae_1)$, $a_3\in \mathbb{R}\backslash \{0\}$.
Since ${\bf k}$ is a $3$-dimensional subalgebra of the $5$-dimensional abelian ideal ${\bf v}=\langle f_{1},e_{1},e_{2},e_{3},e_4\rangle$ it has the form
${\bf k}=\langle e_{2}+b_1 f_1+a_1 e_1, e_{3}+b_2f_1+ a_2e_1, e_{4}+a_3(f_1+ae_{1})\rangle$, $a, a_i, b_i \in \mathbb{R}$, $i=1,2,3$, $aa_3\neq0$. Using the automorphism $\phi(f_{1})=f_{1}$, $\phi(e_{i})=e_{i}$, $i=1,5$,
$\phi(e_2)=e_2-a_1e_1$, $\phi(e_3)=e_3-a_1e_2-a_2e_1$ and $\phi(e_{4})=a_{3}e_{4}$, the subalgebra ${\bf k}$ reduces to ${\bf k}_{6,3}=\langle e_{2}+b_1f_1,e_{3}+b_2f_1,e_{4}+f_{1}+ae_1\rangle$.
This proves the assertion.
\end{proof}
Using  (\cite{strugar}, \S4) we obtain:

\begin{lemma}\label{Liegroups2}
The linear representation of the simply connected Lie group $G_i$ and its subgroup $K_i$ of the Lie algebra ${\bf g}_i$ and its subalgebra ${\bf k}_i$, $i=1,...,8$, is given by the multiplication: 
For $i=1$
\begin{equation} g(x_1,x_2,x_3,x_4,x_5,x_6) g(y_1,y_2,y_3,y_4,y_5,y_6)= \nonumber \end{equation}
\begin{equation} g(x_1+y_1e^{ax_4},x_2+(y_2+x_4y_3)e^{x_4},x_3+y_3e^{x_4},x_4+y_4,x_5+y_5,x_6+y_6), a\neq 0,\nonumber \end{equation}
\begin{equation} K_{1}=\{g(u_1,u_2,u_3,0,u_1+u_2,0); u_i \in \mathbb R, i=1,2,3 \}, \nonumber \end{equation}

For $i=2$
\begin{equation}g(x_1,x_2,x_3,x_4,x_5,x_6) g(y_1,y_2,y_3,y_4,y_5,y_6)=g(x_1+(y_1+x_4y_2+\frac{1}{2}x_4^{2}y_3)e^{x_4}, \nonumber \end{equation}
\begin{equation} x_2+(y_2+x_4y_3)e^{x_4},x_3+y_3e^{x_4},x_4+y_4,x_5+y_5,x_6+y_6), \nonumber \end{equation}
\begin{equation} K_{2}=\{g(u_1,u_2,u_3,0,u_1+a_2u_2+a_3u_3,0);u_i \in \mathbb R, i=1,2,3 \}, a_2,a_3 \in \mathbb R,\nonumber \end{equation}

For $i=3$
\begin{equation} g(x_1,x_2,x_3,x_4,x_5,x_6) g(y_1,y_2,y_3,y_4,y_5,y_6)= \nonumber \end{equation}
\begin{equation} g(x_1+y_1e^{x_4},x_2+y_2e^{ax_4},x_3+y_3e^{bx_4},x_4+y_4,x_5+y_5,x_6+y_6), -1 \le a \le b \le 1, \nonumber \end{equation}
\begin{equation} ab\neq 0, K_{3}=\{g(u_1,u_2,u_3,0,u_1+u_2+u_3,0); u_i \in \mathbb R, i=1,2,3 \}, \nonumber \end{equation}

For $i=4$
\begin{equation}g(x_1,x_2,x_3,x_4,x_5,x_6)g(y_1,y_2,y_3,y_4,y_5,y_6)=g(x_1+y_1e^{ax_4}, \nonumber \end{equation}
\begin{equation} x_2+(y_2cos (x_4)+y_3sin (x_4))e^{bx_4},x_3+(y_3cos (x_4)-y_2sin (x_4))e^{bx_4}, \nonumber \end{equation}
\begin{equation}x_4+y_4,x_5+y_5,x_6+y_6), a\neq 0, b\geq 0, \nonumber \end{equation}
\begin{equation} K_{4,1}=\{g(u_1,u_2,u_3,0,u_1+u_2+a_3u_3,0);u_i \in \mathbb R, i=1,2,3 \}, a_3 \in \mathbb R,\nonumber \end{equation}
\begin{equation} K_{4,2}=\{g(u_1,u_2,u_3,0,u_1+u_3,0); u_i \in \mathbb R, i=1,2,3 \}, \nonumber \end{equation}

For $i=5$
\begin{equation} g(x_1,x_2,x_3,x_4,x_5,x_6) g(y_1,y_2,y_3,y_4,y_5,y_6)= \nonumber \end{equation}
\begin{equation} g(x_1+y_1+x_5y_2,x_2+y_2,x_3+y_3e^{x_5},x_4+y_4e^{cx_5},x_5+y_5,x_6+y_6), 0 < |c| \le 1,
\nonumber \end{equation}
\begin{equation} K_{5,\epsilon}=\{g(u_2+u_3,u_1,u_2,u_3,0,\epsilon u_1); u_i \in \mathbb R, i=1,2,3 \}, \epsilon=0,1, \nonumber \end{equation}

For $i=6$
\begin{equation} g(x_1,x_2,x_3,x_4,x_5,x_6) g(y_1,y_2,y_3,y_4,y_5,y_6)= \nonumber \end{equation}
\begin{equation} g(x_1+y_1+x_5y_2+\frac{1}{2}x_{5}^2y_3,x_2+y_2+x_5y_3,x_3+y_3,x_4+y_4e^{x_5},x_5+y_5,x_6+y_6), \nonumber \end{equation}
\begin{equation} K_{6,\epsilon}=\{g(u_3,u_1,u_2,u_3,0,\epsilon u_2); u_i \in \mathbb R, i=1,2,3 \}, \epsilon=0,1, \nonumber \end{equation}
\begin{equation} K_{6,3}=\{g(au_3,u_1,u_2,u_3,0,b_1u_1+b_2u_2+u_3); u_i \in \mathbb R, i=1,2,3\}, b_1,b_2\in \mathbb R, a\neq 0 ,\nonumber \end{equation}

For $i=7$
\begin{equation} g(x_1,x_2,x_3,x_4,x_5,x_6) g(y_1,y_2,y_3,y_4,y_5,y_6)=g(x_1+y_1+x_2y_5,x_2+y_2, \nonumber \end{equation}
\begin{equation} x_3+(y_3cos (x_5)-y_4sin (x_5))e^{px_5},x_4+(y_4cos (x_5)+y_3sin (x_5))e^{px_5}, \nonumber \end{equation} 
\begin{equation} x_5+y_5,x_6+y_6), p \neq 0, \nonumber \end{equation}
\begin{equation} K_{7,\epsilon}=\{g(u_2+a_3u_3,u_1,u_2,u_3,0,\epsilon u_1); u_i \in \mathbb R, i=1,2,3 \}, \epsilon=0,1, a_3\in \mathbb R, \nonumber \end{equation}
\begin{equation} K_{7,\delta}=\{g(u_3,u_1,u_2,u_3,0,\varepsilon u_1); u_i \in \mathbb R, i=1,2,3 \}, \delta=0,1, \nonumber \end{equation}

For $i=8$
\begin{equation} g(x_1,x_2,x_3,x_4,x_5,x_6) g(y_1,y_2,y_3,y_4,y_5,y_6)= \nonumber \end{equation}
\begin{equation} g(x_1+(y_1+y_2x_5)e^{x_5},x_2+y_2e^{x_5},x_3+y_3+x_5y_4,x_4+y_4,x_5+y_5,x_6+y_6),
\nonumber \end{equation}
\begin{equation}K_{8,\epsilon}=\{g(u_1,u_2,u_1,u_3,0,\epsilon u_3); u_i \in \mathbb R, i=1,2,3 \}, \epsilon=0,1, \nonumber \end{equation}
\end{lemma}

\begin{proposition} \label{excluded2}
There does not exist any $3$-dimensional connected topological proper loop $L$ having ${\bf g}_6$ as the Lie algebra of its multiplication group and the Lie algebra ${\bf k}_{6,3}$ as the Lie algebra of its inner mapping group.
\end{proposition}
\begin{proof}
We may assume that $L$ is simply connected and hence it is homeomorphic to $\mathbb R^3$ (cf. Lemma \ref{simplyconnected}). We show that the Lie group $G_6$ does not allow continuous left transversals $S$ and $T$ to the subgroup $K_{6,3}$ such that for all $s\in S$
and $t\in T$ one has $s^{-1}t^{-1}st \in K_{6,3}$ and the set $S\cup T$ generates $G_6$.

Two arbitrary left transversals to the group $K_{6,3}$ in $G_6$ are:
\begin{equation} S =\{g(u,h_1(u,v,w),h_2(u,v,w),h_3(u,v,w),v,w); u,v,w \in \mathbb R \},  \nonumber \end{equation}
\begin{equation} T =\{g(k,g_1(k,l,m),g_2(k,l,m),g_3(k,l,m),l,m); k,l,m \in \mathbb R \},  \nonumber \end{equation}
where $h_i(u,v,w): \mathbb R^3 \to \mathbb R$ and $g_i(k,l,m): \mathbb R^3 \to \mathbb R$, $i=1,2,3$, are
continuous functions with $h_i(0,0,0)= g_i(0,0,0)= 0$.
The products $s^{-1}t^{-1}st$, $s \in S$, $t \in T$, are elements of $K_{6,3}$ if and only if the equations
$$a(g_3(k,l,m)e^{-l}(1-e^{-v})-h_3(u,v,w)e^{-v}(1-e^{-l}))=vg_1(k,l,m)-vlg_2(k,l,m)-$$
\begin{equation} \label{g6one} lh_1(u,v,w)+lvh_2(u,v,w)+\frac{1}{2}l^{2}h_2(u,v,w)-\frac{1}{2}v^{2}g_2(k,l,m),
\end{equation}
\begin{equation} \label{g6one2} g_3(k,l,m)e^{-l}(1-e^{-v})-h_3(u,v,w)e^{-v}(1-e^{-l})=b_1lh_2(u,v,w)-b_1vg_2(k,l,m) \end{equation}
are satisfied for all $k,l,m,u,v,w\in \mathbb R$.
Applying equation (\ref{g6one2}) equation (\ref{g6one}) becomes simplified to
\begin{equation} \label{g6three} vg_1(k,l,m)+a b_1 vg_2(k,l,m)-vl g_2(k,l,m)-\frac{1}{2} v^2 g_2(k,l,m)= \nonumber \end{equation}
\begin{equation}  lh_1(u,v,w)+ a b_1 l h_2(u,v,w)- lvh_2(u,v,w)-\frac{1}{2}l^{2}h_2(u,v,w). \end{equation}
Using the new functions $g'_1(k,l,m)=g_1(k,l,m)+ a b_1 g_2(k,l,m)-l g_2(k,l,m)$,
$h'_1(u,v,w)=h_1(u,v,w)+ a b_1 h_2(u,v,w)- vh_2(u,v,w)$ equation (\ref{g6three}) reduces to
\begin{equation} \label{g6four} v g'_1(k,l,m)- \frac{1}{2} v^2 g_2(k,l,m)=l h'_1(u,v,w) -\frac{1}{2}l^{2}h_2(u,v,w).
\end{equation}
Equation (\ref{g6four}) holds precisely if the functions $g'_1(k,l,m)$ and $g_2(k,l,m)$, respectively $h'_1(u,v,w)$ and $h_2(u,v,w)$ are polynomials of $l$, respectively of $v$ with order at most $2$. Using this, equation
(\ref{g6one2}) is satisfied if and only if its left hand side and its right hand side are $0$. This holds precisely if one has
$g_3(l)=c(e^l-1)$ and $h_3(v)=c(e^v-1)$, where $c$ is a real constant. In this case the set $S \cup T$ does not generate the group $G_6$. Hence by Proposition \ref{kepka} the group $G_6$ and the subgroup $K_{6,3}$ are not the multiplication group and the inner mapping group of $L$. This proves the assertion.
\end{proof}

\begin{theorem} \label{multiplication2}
Let $L$ be a connected simply connected topological proper loop of dimension $3$ such that its multiplication group is a
$6$-dimensional solvable non-nilpotent decomposable Lie group having $2$-dimensional centre. Then the pairs of Lie groups
$(G_i,K_i)$, $i=1, \cdots, 8$, given in Lemma \ref{Liegroups2} are the multiplication groups $Mult(L)$ and the inner mapping
groups $Inn(L)$ of $L$ with the only exception $(G_6, K_{6,3})$.
\end{theorem}
\begin{proof}
The pairs
$(G_i,K_i)$, $i=1, \cdots, 8$, in Lemma \ref{Liegroups2} can occur as the group $Mult(L)$ and the subgroup
$Inn(L)$ of $L$. According to Proposition \ref{excluded2} the pair $(G_6, K_{6,3})$ is excluded. In all other cases we give  continuous left transversals $A_i$, $B_i$ to the subgroup $K_i$, $i=1, \cdots, 8$, which fulfill the requirements of Proposition \ref{kepka}.

Appropriate $K_1$-connected left transversals in the group $G_1^a$ are:
for $a<-1$ and for $a>1$ the sets
\begin{equation} A_{1,1}=\{g(e^{au}(e^{-u}-1),e^{u}(1-e^{-au})+u^{2},u,u,v,w); u, v, w \in \mathbb R\},
\nonumber \end{equation}
\begin{equation} B_{1,1}=\{g(e^{ak}(1-e^{-k}),k^{2}-e^{k}(1-e^{-ak}),k,k,l,m); k, l, m \in \mathbb R\},
\nonumber \end{equation}
for $0<a<1$ and for $-1<a<0$ the sets
\begin{equation} A_{1,2}=\{g(-ue^{au-u},1-e^{u}+ue^{u}(1-e^{-au}),e^{u}-e^{-au+u},u,v,w); u, v, w
\in \mathbb R\},  \nonumber \end{equation}
\begin{equation} B_{1,2}=\{g(ke^{ak-k},1-e^{k}+ke^{k}(e^{-ak}-1),e^{-ak+k}-e^{k},k,l,m); k, l, m
\in \mathbb R\},  \nonumber \end{equation}
for $a=1$ the sets
\begin{equation} A_{1,3}=\{g(w,e^{u}-1-w+u^{2},u,u,v,w); u, v, w \in \mathbb R\},  \nonumber \end{equation}
\begin{equation} B_{1,3}=\{g(l^{2},e^{k}-1-l^{2}+k^{2},k,k,l,m); k, l, m \in \mathbb R\}, \nonumber \end{equation}
for $a=-1$ the sets
\begin{equation} A_{1,4}=\{g(ue^{-2u},e^{u}-1-ue^{u}+ue^{2u},e^{2u}-e^{u},u,v,w); u, v, w \in \mathbb R\},
\nonumber \end{equation}
\begin{equation} B_{1,4}=\{g(-ke^{-2k},e^{k}-1+ke^{k}-ke^{2k},e^{k}-e^{2k},k,l,m); k, l, m \in \mathbb R\}.
\nonumber \end{equation}
Appropriate $K_2$-connected left transversals in $G_2$ are the sets
\begin{equation} A_2=\{g(e^{u}-1-u^{3}+\frac{3}{2}a_2u^{2}+u(a_3-a_2^{2}),a_2u-\frac{3}{2}u^{2},-u,u,v,w); u,v,w \in \mathbb R\}, \nonumber \end{equation}
\begin{equation} B_2=\{g(e^{k}-1+k^{3}-\frac{3}{2}k^{2}a_2+k(a_2^{2}-a_3),\frac{3}{2}k^{2}-a_2k,k,k,l,m); k, l, m \in \mathbb R\}, \nonumber \end{equation}
$a_2, a_3\in \mathbb{R}$.
Appropriate $K_3$-connected left transversals in $G_3^{a,b}$ are: for $-1\leq a=b\leq 1$ the sets
\begin{equation} A_{3,1}=\{g(e^{u}(e^{-au}-1),e^{au}(1-e^{-u})-w,w,u,v,w); u,v,w \in \mathbb R\}, \nonumber \end{equation}
\begin{equation} B_{3,1}=\{g(e^{k}(1-e^{-ak}),e^{ak}(e^{-k}-1)-m,m,k,l,m); k,l,m \in \mathbb R\}, \nonumber \end{equation}
for $-1\leq a<b\leq 1$ the sets
\begin{equation} A_{3,2}=\{g(e^{u-au}-e^{u-bu},e^{au}-e^{au-u},e^{bu-u}-e^{bu},u,v,w); u, v, w \in \mathbb R\}, \nonumber \end{equation}
\begin{equation} B_{3,2}=\{g(e^{k-bk}-e^{k-ak},e^{ak-k}-e^{ak},e^{bk}-e^{bk-k},k,l,m); k, l, m \in \mathbb R\},  \nonumber \end{equation}
where $ab\neq 0$.
An appropriate $K_{4,1}$-connected left transversal in $G_4^{a,b}$ is the set
\begin{equation} A_{4,1}=B_{4,1}=\{g(e^{au-bu}sin(u), \nonumber \end{equation}
\begin{equation} \frac{1}{a_{3}^{2}+1}((e^{bu}-e^{bu-au})(sin(u)-a_3cos(u))+(e^{bu}-cos(u))(cos(u)+a_3sin(u))),\nonumber \end{equation}
\begin{equation}\frac{1}{a_{3}^{2}+1}((e^{bu}-e^{bu-au})(a_3sin(u)+cos(u))+(e^{bu}-cos(u))(a_3cos(u)-sin(u))),\nonumber \end{equation}
\begin{equation} u,v,w); u, v, w \in \mathbb R\}, a_3 \in \mathbb R. \nonumber \end{equation}
An appropriate $K_{4,2}$-connected left transversal in $G_4^{a,b}$ is the set
\begin{equation} A_{4,2}=B_{4,2}=\{g(e^{au-bu}sin(u),(e^{bu-au}-e^{bu})cos(u)+sin(u)(e^{bu}-cos(u)),\nonumber \end{equation}
\begin{equation} cos(u)(e^{bu}-cos(u))-(e^{bu-au}-e^{bu})sin(u),u,v,w); u, v, w \in \mathbb R\},  \nonumber \end{equation}
$a\neq 0$, $b\geq 0$.
Appropriate $K_{5,\epsilon}$-connected left transversals in $G_5^c$ with $\epsilon=0,1$ are: for $c=1$ the sets
\begin{equation} A_5^{c=1}=\{g(u,1-e^{-v},u,ve^{v}-u,v,w); u, v, w \in \mathbb R\},  \nonumber \end{equation}
\begin{equation} B_5^{c=1}=\{g(k,e^{-l}-1,k,-le^{l}-k,l,m); k, l, m \in \mathbb R\},  \nonumber \end{equation}
for $c\neq 1$ the sets
\begin{equation} A_5^{c \neq 1}=\{g(u,e^{-v}-e^{-cv}, -ve^v,ve^{cv},v,w); u, v, w \in \mathbb R\},  \nonumber \end{equation}
\begin{equation} B_5^{c \neq 1}=\{g(k,e^{-cl}-e^{-l},le^{l},-le^{cl},l,m); k, l, m \in \mathbb R\}.  \nonumber \end{equation}
Appropriate $K_{6,\epsilon}$-connected left transversals in $G_6$, where $\epsilon=0,1$, are the sets
\begin{equation} A_6=\{g(u,1-v^{2}-e^{-v},-v,\frac{1}{2}v^{2}e^{v},v,w); u,v,w \in \mathbb R\}, \nonumber \end{equation}
\begin{equation} B_6=\{g(k,l+\frac{1}{2}l^{2}-le^{-l},1-e^{-l},-le^{l},l,m); k,l,m \in \mathbb  R\}. \nonumber \end{equation}
An appropriate $K_{7,\epsilon}$-connected left transversal in $G_{7}^{p\neq 0}$ with $\epsilon=0,1$ is the set
\begin{equation} A_7=B_7=\{g(u,e^{-pv}sin(v), \nonumber \end{equation}
\begin{equation} \frac{1}{a_{3}^{2}+1}(e^{pv}v(sin(v)+a_3cos(v))+(e^{pv}-cos(v))(cos(v)-a_3sin(v))),\nonumber \end{equation}
\begin{equation}\frac{1}{a_{3}^{2}+1}(e^{pv}v(a_3sin(v)-cos(v))+(e^{pv}-cos(v))(a_3cos(v)+sin(v))),\nonumber \end{equation}
\begin{equation}v,w); u, v, w \in \mathbb R\},  a_3\in \mathbb R. \nonumber \end{equation}
An appropriate $K_{7,\delta}$-connected left transversal in $G_{7}^{p\neq 0}$, where $\delta=0,1$, is the set
\begin{equation} A_7=B_7=\{g(u,e^{-pv}sin(v),e^{pv}vcos(v)+sin(v)(e^{pv}-cos(v)),\nonumber \end{equation}
\begin{equation} e^{pv}vsin(v)-cos(v)(e^{pv}-cos(v)),v,w); u, v, w \in \mathbb R\}.  \nonumber \end{equation}
Appropriate $K_{8,\epsilon}$-connected left transversals in $G_8$ with $\epsilon=0,1$ are the sets
\begin{equation} A_8=\{g(ve^{v}+v^{2},v,u,1-e^{-v},v,w); u,v,w \in \mathbb R\}, \nonumber \end{equation}
\begin{equation} B_8=\{ g(l^{2}-le^{l},l,k,e^{-l}-1,l,m); k,l,m \in \mathbb  R\}. \nonumber \end{equation}
Hence the assertion follows from Proposition \ref{kepka}.
\end{proof}

\begin{corollary} \label{cor1}
Each $6$-dimensional solvable decomposable Lie group which is the group $Mult(L)$ of a $3$-dimensional connected topological loop has $1$- or $2$-dimensional centre and $3$-dimensional commutator subgroup.
\end{corollary}
\begin{proof} If $L$ has $1$-dimensional centre, then the assertion follows from Proposition \ref{1dimcentre}. If $L$ has $2$-dimensional centre, then Proposition \ref{2dimensionalcentre} yields the assertion. \end{proof}

\begin{corollary} \label{cor2}
Each $3$-dimensional connected topological proper loop $L$ having a solvable Lie group of dimension $6$ as the group
$Mult(L)$ of $L$ has $1$- or $2$-dimensional centre and $2$- or $3$-dimensional commutator subgroup.
\end{corollary}
\begin{proof}
If $L$ has a $6$-dimensional solvable indecomposable Lie group as its multiplication group, then the assertion is proved in Corollary 3.4 in \cite{ffaa}. If $L$ has a $6$-dimensional solvable decomposable Lie group as its multiplication group, then
Corollary \ref{cor1} gives the assertion.
\end{proof}

\end{document}